\documentclass{gtart_h}

\def\ifplaintex{\expandafter\ifx\csname documentclass\endcsname\relax}

\def\gtp{{\mathsurround=0pt\it $\cal G\mskip-2mu$eometry \&\ 
$\cal T\!\!$opology $\cal P\!$ublications}}  

\def\recd{{\small Received:\qua\receiveddate\ifx\reviseddate\relax
\else\qquad Revised:\qua\reviseddate\fi\par}} 

\def\lognumber#1{\def\thelognumber{#1}}
\def\volumenumber#1{\def\thevolumenumber{#1}}
\def\volumeyear#1{\def\thevolumeyear{#1}}
\def\papernumber#1{\def\thepapernumber{#1}}
\def\pagenumbers#1#2{\def\startpage{#1}\def\finishpage{#2}}
\def\published#1{\def\publishdate{#1}}

\def\received#1{\def\receiveddate{#1}}
\def\revised#1{\def\reviseddate{#1}}
\def\accepted#1{\def\accepteddate{#1}}
\def\asciititle#1{\def\theasciititle{#1}}

\def\asciiauthors#1{\def\theasciiauthors{#1}}

\def\coverauthors#1{\def\thecoverauthors{#1}}
\long\def\asciiabstract#1{\long\def\theasciiabstract{#1}}
\def\asciikeywords#1{\def\theasciikeywords{#1}}

\let\\\par\let\thelognumber\relax\let\thevolumenumber\relax
\let\thepapernumber\relax\let\thevolumeyear\relax\let\startpage\relax
\let\finishpage\relax\let\publishdate\relax\let\receiveddate\relax
\let\reviseddate\relax\let\accepteddate\relax\let\theasciititle\relax
\let\theasciiauthors\relax
\let\theasciiabstract\relax\let\theasciikeywords\relax

\let\thecoverauthors\relax\let\theasciiemail\relax

\ifplaintex
\font\logobig=cmssbx10 scaled 3836
\font\logomed=cmssbx10 scaled 2557
\else
\font\logobig=cmssbx10 scaled 4200
\font\logomed=cmssbx10 scaled 2800
\fi

\long\def\makeagttitle{   
\count0=\startpage
\agt\hfill      
\hbox to 45truept{\vbox to 0pt{\vglue -13truept{\logomed A\kern -.37em{\logobig 
T}\kern -.38em G}\vss}\hss}
\break
{\small Volume \thevolumenumber\ (\thevolumeyear)
\startpage--\finishpage\nl
Published: \publishdate}

\vglue .25truein

{\parskip=0pt\leftskip 0pt plus
1fil\def\\{\par\smallskip}{\Large\bf\thetitle}\par\medskip} \vglue
0.05truein

{\parskip=0pt\leftskip 0pt plus 1fil\def\\{\par}{\sc\theauthors}
\par\medskip}%
 
\vglue 0.03truein 

{\small\leftskip 25truept\rightskip 25truept{\bf Abstract}\stdspace\theabstract

{\bf AMS Classification}\stdspace\theprimaryclass
\ifx\thesecondaryclass\relax\else; \thesecondaryclass\fi\par
{\bf Keywords}\stdspace \thekeywords\par}\vglue 7truept

}   

\ifplaintex
\font\phead=cmsl9 scaled 950
\font\pnum=cmbx10 scaled 913
\font\pfoot=cmsl9 scaled 950
\headline{\vbox to 0pt{\vskip -4.5mm\line{\small\phead\ifnum
\count0=\startpage ISSN 1472-2739 (on-line) 1472-2747 (printed)
\hfill {\pnum\folio}\else\ifodd\count0\def\\{ }%
\ifx\theshorttitle\relax\thetitle\else\theshorttitle\fi\hfill{\pnum\folio}
\else\def\\{ and }{\pnum\folio}\hfill\ifx\theshortauthors\relax\theauthors
\else\theshortauthors\fi\fi\fi}\vss}}
\footline{\vbox to 0pt{\vglue 0mm\line{\small\pfoot\ifnum\count0=\startpage
\copyright\ \gtp\hfill\else
\agt, Volume \thevolumenumber\ (\thevolumeyear)\hfill\fi}\vss}}
\else
\font=cmsl9 scaled 1050
\font\lnum=cmbx10 
\font=cmsl9 scaled 1050
\makeatletter
\def\@oddhead{{\small\lhead\ifnum\count0=\startpage ISSN 1472-2739 
(on-line) 1472-2747 (printed)\hfill {\lnum\number\count0}\else\ifodd\count0
\def\\{ }\ifx\theshorttitle\relax \thetitle \else\theshorttitle\fi\hfill
{\lnum\number\count0}\else\def\\{ and }{\lnum\number\count0}
\hfill\ifx\theshortauthors\relax 
\theauthors\else\theshortauthors\fi\fi\fi}}\def\@evenhead{\@oddhead}
\def\@oddfoot{\small\lfoot\ifnum\count0=\startpage\copyright\ \gtp\hfill\else
\agt, Volume \thevolumenumber\ (\thevolumeyear)\hfill\fi}
\def\@evenfoot{\@oddfoot}
\makeatother
\fi
\let\maketitlepage\makeagttitle
\let\makeshorttitle\maketitlepage
\let\maketitle\maketitlepage

\newwrite\gtoutfile
\long\gdef\makeheadfile{  
{\def\\{, }\def\s{ }
\immediate\openout\gtoutfile head.xxx
\immediate\write\gtoutfile{Proxy-for: \ifx\theasciiauthors\relax
\theauthors\else\theasciiauthors\fi\s<\ifx\theasciiemail\relax\theemail\else\theasciiemail\fi>}
\immediate\write\gtoutfile{\noexpand\\}
\immediate\write\gtoutfile{Authors: \ifx\theasciiauthors\relax
\theauthors\else\theasciiauthors\fi}
{\def\\{ }\immediate\write\gtoutfile{Title: \ifx\theasciititle\relax
\thetitle\else\theasciititle\fi}}
\immediate\write\gtoutfile{Subj-class: GT or SG, GR etc}
\immediate\write\gtoutfile{MSC-class: \theprimaryclass\ifx\thesecondaryclass\relax\else, \thesecondaryclass\fi}
\immediate\write\gtoutfile{Journal-ref: Algebr. Geom. Topol. \thevolumenumber\s
(\thevolumeyear) \startpage-\finishpage}
\immediate\write\gtoutfile{Comments: Published by Algebraic and
Geometric Topology at}
\immediate\write\gtoutfile{\s\s\s  http://www.maths.warwick.ac.uk/agt/AGTVol\thevolumenumber/agt-\thevolumenumber-\thepapernumber.abs.html}
\immediate\write\gtoutfile{\noexpand\\}
\immediate\write\gtoutfile{}
\ifx\theasciiabstract\relax
\immediate\write\gtoutfile{\theabstract}\else
\immediate\write\gtoutfile{\theasciiabstract}\fi
\immediate\write\gtoutfile{}
\immediate\write\gtoutfile{\noexpand\\}
\immediate\write\gtoutfile{}
\immediate\closeout\gtoutfile}}  

\def\maketitlepage{\makeagttitle\makeheadfile}
\let\makeshorttitle\maketitlepage
\let\maketitle\maketitlepage

\lognumber{65}
\volumenumber{5}
\volumeyear{2005}
\papernumber{65}
\pagenumbers{1589}{1635}
\received{15 June 2005} 
\revised{9 November 2005}  
\accepted{15 November 2005}
\published{24 November 2005}

\usepackage{amssymb, amsmath, mathrsfs}
\usepackage[all]{xy}

\newcommand{\Z}{{\mathbb Z}}

\newcommand{\F}{{\mathbb F}}

\def\:{\co}
\def\.{\cdot}
\def\<{\left\langle}
\def\>{\right\rangle}
\def\({\left(}
\def\){\right)}
\def\epsilon{\varepsilon}
\def\phi{\varphi}
\def\subset{\subseteq}

\def\leq{\leqslant}
\def\geq{\geqslant}

\def\lra{\longrightarrow}
\def\Lra{\Longrightarrow}

\def\Mod{\mathsf{Mod}}
\def\Comod{\mathsf{Comod}}

\def\dt{\mathscr{D}_T}
\def\dr{\mathscr{D}_R}
\def\ds{\mathscr{D}_{\mathbb S}}
\def\ms{\mathscr{M}_{\mathbb S}}

\def\xra{\xrightarrow}

\def\calc{\mathcal C}

\def\cale{\mathcal E}
\def\calf{\mathcal F}
\def\calj{\mathcal J}

\def\wt{\widetilde}
\def\wh{\widehat}
\def\ss{\mathbb S}

\newtheorem{thm}{Theorem}\newtheorem{lem}[thm]{Lemma}
\newtheorem{prop}[thm]{Proposition}\newtheorem{cor}[thm]{Corollary}

\numberwithin{equation}{section} \numberwithin{thm}{section}

\theoremstyle{remark}
\newtheorem{rem}[thm]{Remark}
\newtheorem{exmp}[thm]{Example}

\theoremstyle{definition}
\newtheorem{defn}[thm]{Definition}
\newtheorem{notation}[thm]{Notation}

\DeclareMathOperator{\Hom}{Hom} 
 \DeclareMathOperator{\gr}{gr}
\DeclareMathOperator{\Ext}{Ext} \DeclareMathOperator{\Tor}{Tor}

\DeclareMathOperator{\coker}{coker}
\DeclareMathOperator{\hocolim}{hocolim}\DeclareMathOperator{\holim}{holim}
\DeclareMathOperator{\colim}{colim}
\DeclareMathOperator{\Coext}{Coext}
\DeclareMathOperator{\Cohom}{Cohom}
 \DeclareMathOperator{\Sym}{Sym}
\DeclareMathOperator{\RHom}{{\mathbf R}Hom}
\DeclareMathOperator{\hatotimes}{\widehat\otimes}
\DeclareMathOperator{\derotimes}{\otimes^{\mathbf L}}
\DeclareMathOperator{\GL}{GL}

\newcommand{\ie}{i.e.}
\newcommand{\eg}{e.g.}

\begin{document}
\title{$I$--adic towers in topology}
\asciititle{I-adic towers in topology}

\author{Samuel W\"uthrich}
\coverauthors{Samuel W\noexpand\"uthrich}
\asciiauthors{Samuel Wuethrich}

\address{Department of Pure Mathematics, University of 
Sheffield\\Hicks Building, Hounsfield Road, Sheffield S3 7RH, UK}

\email{S.Wuethrich@sheffield.ac.uk}

\begin{abstract}
A large variety of cohomology theories is derived from complex
cobordism $MU^*(-)$ by localizing with respect to certain elements
or by killing regular sequences in $MU_*$. We study the
relationship between certain pairs of such theories which differ
by a regular sequence, by constructing topological analogues of
algebraic $I$--adic towers. These give rise to Higher Bockstein
spectral sequences, which turn out to be Adams spectral sequences
in an appropriate sense. Particular attention is paid to the case
of completed Johnson--Wilson theory $\wh{E(n)}$ and
Morava $K$-theory $K(n)$ for a given prime $p$.
\end{abstract}

\asciiabstract{%
A large variety of cohomology theories is derived from complex
cobordism MU^*(-) by localizing with respect to certain elements or by
killing regular sequences in MU_*. We study the relationship between
certain pairs of such theories which differ by a regular sequence, by
constructing topological analogues of algebraic I-adic
towers. These give rise to Higher Bockstein spectral sequences, which
turn out to be Adams spectral sequences in an appropriate
sense. Particular attention is paid to the case of completed
Johnson--Wilson theory E(n)-hat and Morava K-theory K(n) for a
given prime p.}

\primaryclass{55P42, 55P43, 55T15}
\secondaryclass{55U20, 55P60, 55N22}

\keywords{Structured ring spectra, Adams resolution, Adams
spectral sequence, Bockstein operation, complex cobordism, Morava
$K$-theory, Bousfield localization, stable homotopy theory.}
\asciikeywords{Structured ring spectra, Adams resolution, Adams
spectral sequence, Bockstein operation, complex cobordism, Morava
K-theory, Bousfield localization, stable homotopy theory.}

\maketitle

\section*{Introduction}

Complex cobordism $MU^*(-)$ is arguably one of the most powerful
cohomology theories. It describes stable phenomena with high
accuracy. A famous piece of evidence for such a vague statement
are the nilpotence and periodicity theorems \cite{ravenelorange}.
These lead to the chromatic filtration of the stable homotopy
groups of spheres. Due to its strong geometric background, the
spectrum $MU$ representing complex cobordism is multiplicative in
a highly structured way. Namely, $MU$ can be realized as an
$E_\infty$ ring spectrum \cite{mayeinfty} or equivalently, in the
framework of \cite{ekmm}, as a commutative $\ss$--algebra. As a
drawback to the large amount of information that $MU^*(X)$ carries
about a given space $X$, it is rather difficult to explicitly
determine $MU^*(X)$ in general. The approach to confront this
problem is to construct new cohomology theories by formally
manipulating the theory $MU^*(-)$. Important possibilities are to
realize regular quotients and localizations of the coefficient
ring $MU_*$ as the coefficient rings of new theories.

In modern geometric categories like the ones constructed in
\cite{ekmm}, it is possible to perform such constructions in a
transparent way. Namely, thanks to the fact that $MU$ is a
commutative $\ss$--algebra, we may work in the monoidal and
triangulated category $\mathscr{D}_{MU}$ of $MU$--modules. The
monoidal structure is provided by the smash product $\wedge_{MU}$
over $MU$. Furthermore, it is possible to construct the new
spectra as ring objects (monoids) in $\mathscr{D}_{MU}$, so-called
$MU$--ring spectra \cite{strickland}.

A natural question arises: What is the relationship among theories
constructed in such a fashion? In this paper, we study the
connections between two $MU$--ring spectra $T$ and $F$, where $F$
is obtained from $T$ by killing an ideal $I$ of $T_*$, which is
generated by a regular sequence $S$. In order to avoid
technicalities, we don't make this more precise now. More
generally, we consider ring spectra $T$ and $F$ over an arbitrary
commutative $\ss$--algebra $R$. An example for $R=MU$ is completed
Johnson--Wilson theory $\wh{E(n)}$ as $T$ and Morava $K$-theory
$K(n)$ as $F$, for a specified prime $p$. Given a pair $T$ and
$F$, we can ask whether there is a way to compute $T^*(X)$ from
$F^*(X)$. If the sequence $S$ is finite there is a standard
procedure to do this. Namely, we can work our way through a
sequence of Bockstein spectral sequences. For the special case
$\wh{E(n)}$ and $K(n)$ for an odd prime $p$, Baker and W\"urgler
describe in \cite{bakainfty, bw} a different approach. They
realize the $I$--adic tower under $T_*$
\begin{equation}\label{iadicintro}
\begin{array}{c}
\xymatrix
{\quad T_*\quad \ar[r] & \quad \cdots \quad
\ar[r] & T_*/I^3 \ar[r] & T_*/I^2 \ar[r] & F_*=T_*/I
\\
&& I^2/I^3 \ar[u] & I/I^2 \ar[u]}
\end{array}
\end{equation}
in topology and derive the so-called Higher Bockstein spectral
sequence with $E_1$--term depending only on $F^*(X)$ and target
$T^*(X)$. Based on the technology from \cite{ekmm}, Baker and
Lazarev \cite{bl} discuss a construction of $I$--adic towers for
the case where the ring spectrum $T$ is a commutative
$\ss$--algebra. They note a close relationship of such towers to
Adams resolutions.

The aim of this paper is to provide an alternative construction of
$I$--adic towers which the author believes is more transparent and
which explains their conceptual nature. This is achieved by making
precise in what sense such towers are Adams resolutions. Besides
the theoretical desirability of understanding the nature of
$I$--adic towers, there are several practical advantages. Perhaps
most importantly, the spectra $T/I^s$ constituting the tower
become much more tractable. One concrete application the author
has in mind is the study of multiplicative structures on these
objects \cite{swmult}. Another benefit is that it becomes possible
to specify clearly what the construction depends on. It turns out
that one can considerably weaken the assumption that the ring
spectrum $T$ be strictly commutative. This will be made more
explicit further below.

The strategy pursued here to develop an appropriate concept of a
topological $I$--adic tower is to analyze the algebraic model in
the derived setting. Namely, suppose that $T_*$ is commutative and
consider the $I$--adic tower over $T_*$ in the derived category
$\mathscr{D}_{T_*}$. It is induced by the diagram
\eqref{iadicintro} of $T_*$--modules and is of the form:
\begin{equation}\label{iadicintro2}
\begin{array}{c}
\xymatrix@!C=.8cm{ \cdots\ \ar[r] & I^3 \ar[r]\ar[d] & I^2
\ar[r]\ar[d] & I \ar[r]\ar[d] & T_*\ar[d]
\\
& I^3/I^4 \ar[lu]|-\circ & I^2/I^3 \ar[lu]|-\circ & I/I^2
\ar[lu]|-\circ & T_*/I=F_* \ar[lu]|-\circ }
\end{array}
\end{equation}
Arrows with a circle denote maps of degree $-1$. The fundamental
observation made in \cite{swshort} is that the sequence
\begin{equation}\label{ressequence}
\begin{array}{c}
\xymatrix{ T_* \ar[r] & F_* \ar[r]|-\circ & I/I^2 \ar[r]|-\circ &
I^2/I^3 \ar[r]|-\circ  & \ \cdots}
\end{array}
\end{equation}
derived from \eqref{iadicintro2}, determines the tower up to
isomorphism. The point is that \eqref{ressequence} is a relative
injective resolution of $T_*$ with respect to $F_*$ (in the sense
of Remark \ref{adamsalgebraic}) whose associated Adams resolution
is given by \eqref{iadicintro2} (Remark
\ref{adamsalgebraicresolution}).

To imitate such a construction in topology, we need to have enough
structure on our objects. We define the notion of a regular triple
$(R, T, F)$ of ring spectra (Definition \ref{regulartriple}),
which is meant to incorporate a rather natural and general setup
in which it makes sense to consider the problem of constructing
$I$--adic towers. In fact, we show that it is possible to
construct such towers for any given regular triple. In the
following, we motivate the ingredients which go into the
definition. As mentioned before, $R$ is assumed to be a
commutative $\ss$--algebra. It plays the role of the ground ring,
in the sense that we work in the category of $R$--modules $\dr$,
equipped with the smash product $\wedge_R$. We suppose that $T$ is
a given commutative $R$--ring spectrum and that $F$ is the
quotient of $T$ by a $T_*$-regular sequence $S=(x_0, x_1, \ldots)$
in $R_*$. By this, we mean that $F$ is of the form $T\wedge_R L$
with $L=R/x_0\wedge_R R/x_1 \wedge_R \cdots$. Here $R/x_i$ denotes
the cofibre of the multiplication map $x_i$, regarded as a
(graded) endomorphism of $R$. Then the homotopy groups of $F$
realize the quotient $T_*/I$, where $I$ is the ideal generated by
the image of $S$ under $(\eta_T)_*\: R_*\to T_*$, induced by the
unit map $\eta_T\: R\to T$ of $T$. To obtain a manageable notion
of relative injective resolutions in $\dr$ with respect to $F$, we
need $F$ to be an $R$--ring spectrum (see Remark
\ref{adamsmonoid}). Rather than assuming this to be part of the
data, we specify conditions which guarantee that there is an
$R$--ring structure on $L$ and hence one on $F$. Namely, we assume
that the coefficient ring $R_*$ of $R$ is trivial in odd degrees,
that the sequence $S$ is regular on $R_*$ and that it consists of
non-zero divisors. Strickland \cite{strickland} proves that these
conditions are sufficient for our purpose.

In the special case $T=R$, the construction of the topological
analogue of the sequence \eqref{ressequence} is rather
straightforward. We prove in Theorem \ref{iadic} that it is a
relative injective resolution with respect to $F$ (Definition
\ref{defadams}). This allows us to define a topological $I$--adic
tower as the Adams resolution associated to it. We show that its
diagram of homotopy groups realizes \eqref{iadicintro2}. From
there, we may construct a topological $I$--adic tower realizing
\eqref{iadicintro}:
\begin{equation}\label{iadicintro3}
\begin{array}{c}
\xymatrix@!C=1.3cm{\quad T\quad \ar[r] & \quad \cdots \quad
\ar[r] & T/I^3 \ar[r] & T/I^2 \ar[r] & F=T/I\\
&& I^2/I^3 \ar[u] & I/I^2 \ar[u]}
\end{array}
\end{equation}
For $T\neq R$, we can do the construction above for the triple
$(R, R, L)$ and apply $T\wedge_R -$ to it. This is an Adams
resolution, and we prove that its homotopy groups realize
\eqref{iadicintro2}.

If the sequence $S$ is finite, we show that the homotopy limit
$\wh T = \holim_s T/I^s$ is isomorphic to the Bousfield
localization $L_F^R T$ of $T$ with respect to $F$ in $\dr$
(Proposition \ref{localization}). If furthermore $I$ is invariant
in $T_*(T)$, \ie{} if we have $I\cdot T_*(T)=T_*(T)\cdot I$ as
subgroups of the $T_*$--bimodule $T_*(T)$, $\wh T$ is isomorphic
to the Bousfield localization $L_F T$ in the ordinary stable
homotopy category. For $T=E(n)$, we identify $\wh{E(n)}$ as the
spectrum representing completed Johnson--Wilson theory from
\cite{bwart} (Proposition \ref{enhat}) and thus recover the result
$L_{K(n)} E(n) \cong \wh{E(n)}$ proved there.

From \eqref{iadicintro3}, we derive the Higher Bockstein spectral
sequences (Theorem \ref{hbss}):
\begin{align}
\label{cohomhbssintro} E_1^{*,*}  & =  \bigoplus_{s\geq 0}
I^s/I^{s+1} \otimes_{F^*} F^*(X)  \Longrightarrow   \wh T^*(X)
\\
\label{homhbssintro} E^1_{*,*} & = \bigoplus_{s\geq 0} I^s/I^{s+1}
\otimes_{F_*} F_*(X) \Longrightarrow (F^\wedge_R (T\wedge X))_*
\end{align}
The target of \eqref{homhbssintro} are the homotopy groups of the
so-called $F$--nilpotent completion of $T\wedge X$ in $\dr$. We
identify the first differentials as expressions in certain
Bockstein operations. In case that $F$ is commutative and
$F_*^R(L)$ is $F_*$--flat, we identify the $E_2$--terms as certain
$\Coext$--groups over an exterior coalgebra, which arises
naturally as $F_*^R(L)$. For $T=\wh{E(n)}$ and $F=K(n)$ for a
given prime $p$, we show that the target of \eqref{homhbssintro}
is $(L_{K(n)}(\wh{E(n)}\wedge X))_*$ for any $X$ (Proposition
\ref{enhatnilpotent}).

The paper is organized as follows. In Section \ref{notation}, we
recall notation and terminology from \cite{ekmm}. In Section
\ref{kunnethsec}, we discuss duality between homology and
cohomology theories represented by a not necessarily commutative
$R$--ring spectrum $F$ in $\dr$, where $R$ is a commutative
$\ss$--algebra. Furthermore, we record K\"unneth isomorphisms
which are adapted from \cite{boardmanstable} to our situation. In
Section \ref{regulartriples}, we define the notion of a regular
triple, give a list of examples and do some preparatory work. In
Section \ref{algebraiadic}, we collect the facts that we need
about algebraic $I$--adic towers from \cite{swshort, swglasgow}.
In Section \ref{adamsresolutions}, we recall the theory of Adams
resolutions and specify the injective class in $\mathscr{D}_R$
needed in Section \ref{topologytower} to construct $I$--adic
towers. The construction itself is quite direct. Technically more
involved is the verification that the homotopy and $F$--homology
groups are as expected. Once this is done, we identify $\wh T$ as
a Bousfield localization in the sense and under the conditons
mentioned above. We give a direct proof that for $T=E(n)$ and
$I=I_n$, the homotopy limit $\wh T$ represents completed
Johnson--Wilson theory, as defined in \cite{bwart}. In Section
\nolinebreak \ref{sectionhbss}, we discuss the Higher Bockstein
spectral sequences, derived from $I$--adic towers. For the
identification of the $E_2$--terms under the conditions already
mentioned, we decompose the natural $F_*^R(F)$--coaction into
compatible coactions of $F_*^R(T)$ and $F_*^R(L)$. Only the
coaction of the coalgebra $F_*^R(L)$ is relevant for the
$E_2$--terms.  By comparing explicit descriptions of the Bousfield
localization functors $L_{K(n)}^{\wh{E(n)}}$ and $L_{K(n)}$ from
\cite{greenleesmay} and \cite{hs} respectively, we identify the
target of \eqref{homhbssintro} for $T=\wh{E(n)}$ and $F=K(n)$ as
stated above. We end by considering some examples.

This paper forms the main part of my PhD thesis. I am deeply
indebted to my supervisor Alain Jeanneret and my temporary
supervisor Andrew Baker for all their encouragement and support. I
would like to express my gratitude to the following persons: Urs
W\"urgler for proposing the project to analyze Higher Bockstein
spectral sequences; John Greenlees, John Rognes and Neil
Strickland for helpful comments and the referee for pointing out a
subtlety in the definition of regular triples. I am grateful to
the Swiss National Science Foundation for financial support.

\section{Notation and terminology}\label{notation}

We use the framework of \cite{ekmm} in this paper. We start by
recalling some notation and terminology from there.

The category of {\em $\ss$--modules} $\ms$ is a symmetric monoidal
model category with smash product $\wedge_\ss$ and unit the sphere
spectrum $\ss$. Its homotopy category $\mathscr{D}_\ss$, which has
the same objects as $\mathscr M_\ss$, is isomorphic to the
classical stable homotopy category. When regarding an
$\ss$--module as an object of $\mathscr D_\ss$, we call it a {\em
spectrum}. The homotopy groups $\pi_*(E)=[\ss, E]_*$ of a spectrum
$E$ are written as $E_*$. We follow the convention
$(-)^*=(-)_{-*}$ for graded objects. We do not distinguish a map
from its suspensions in the notation.

A monoid $R$ in $\ms$ is called an {\em $\ss$--algebra}. For any
such, we have the model category of (strict) left $R$--modules
$\mathscr{M}_R$. Its associated homotopy category $\mathscr{D}_R$
is a triangulated category. The graded $R^*$--module of morphisms
between two objects $M$ and $N$ is denoted by $[M,N]_R^*$ or
$N^*_R(M)$. In particular, we have $N_* \cong [R, N]_R^{-*}$. The
homological degree of an element $v\in M_*$ is denoted by $|v|$.
We also have the strict and the homotopy categories of right
$R$--modules and $R$--bimodules. When speaking of an {\em
$R$--module}, we mean a left $R$--module $M$, considered as an
object of $\mathscr{D}_R$. When we view $M$ as an object in
$\mathscr M_R$, we refer to it as a {\em strict $R$--module}.

The smash product $\wedge_R$ over $R$ is a bifunctor, which
assigns to a strict right and a strict left $R$--module $M$ and
$N$ respectively an $\ss$--module $M\wedge_R N$. The function
$R$--module functor sends a pair $N$ and $L$ of strict left
$R$--modules to an $\ss$--module $F_R(N,L)$. There is an obvious
adjunction formula between $\wedge_R$ and $F_R$. Both functors
descend to homotopy categories. We have $\pi_{-*}(F_R(N, L)) \cong
[N, L]_R^*$ and write $M_*^R(N)$ for $\pi_*(M\wedge_R N)$. The
functors $L_R^*(-)$ and $M_*^R(-)$ define cohomology and homology
functors on $\mathscr{D}_R$ respectively. If $R$ is clear from the
context, we refer to them as {\em $L$--cohomology \/} and {\em
$M$--homology}.

The category of strict $R$--bimodules is a monoidal model
category; hence its homotopy category is monoidal. If $R$ is a
commutative $\ss$--algebra, the strict categories of left and
right $R$--modules are isomorphic and every strict one-sided
module is a strict bimodule. Hence $\wedge_R$ gives rise to a
monoidal model category structure on $\mathscr{M}_R$ and therefore
to a monoidal structure on $\mathscr{D}_R$. Monoids in
$\mathscr{M}_R$ and $\mathscr{D}_R$ are called {\em
$R$--algebras\/} and {\em $R$--ring spectra\/} respectively.

If $T$ is an $R$--ring spectrum, we may consider algebras over the
monad \cite[VI]{maclane} $T\wedge_R-$ in $\dr$. These are objects
$M$ of $\dr$ with an action map $T\wedge_R M\to M$ in $\dr$,
satisfying the associativity and unit axioms. We refer to such $M$
as {\em $T$--module spectra}.

For a commutative $\ss$--algebra $R$ and an $R$--module $M$,
multiplication by an element $x\in R_*$ defines a map
$\Sigma^{|x|} M\xrightarrow{x} M$ in $\mathscr{D}_R$. We use the
notation
\begin{equation}\label{standardcofibre}
\Sigma^{|x|} M\xrightarrow{\ v\ } M \xrightarrow{\ \rho_x\ } M/x
\xrightarrow{\ \beta_x\ } \Sigma^{|x|+1} M
\end{equation}
from \cite{strickland} for its cofibre and the maps to and from
it. Note that $M/x \cong M\wedge_R R/x$. Furthermore, we define
$\bar\beta_x$ to be the composition
\begin{equation}\label{standardbockstein}
\bar\beta_x = \rho_x\circ\beta_x \: M/x \lra \Sigma^{|x|+1} M/x.
\end{equation}
For a sequence $S=(x_1,x_2, \ldots)$ of elements of $R_*$, we
write
\[
M/S = M\wedge\bigwedge_{j} R/x_j
\]
and call $M/S$ a {\em quotient\/} of $M$. Here, as well as
throughout the paper, $\wedge$ denotes $\wedge_R$, where $R$ is
specified in the context. The symbol $\ast$ denotes a point. A
blank $\otimes$ means $\otimes_{F_*}$, where $F_*$ is the
coefficient ring of some ring spectrum $F$ under consideration.
The category of $F_*$--modules is denoted by $\Mod_{F_*}$.

\section{K\"unneth isomorphisms for $R$--ring spectra}
\label{kunnethsec}

Let $R$ be a commutative $\ss$--algebra and let $F$ be an
$R$--ring spectrum. We can define an exterior product
\begin{equation}\label{kunneth}
\times\: F_*^R(M)\otimes F_*^R(N) \lra F_*^R(M\wedge N)
\end{equation}
for given $R$--modules $M$ and $N$ in the usual way. However, as
$F$ is not assumed to be commutative, we need to take care with
the right action of $F_*$ on $F_*^R(M)$. For $x\in F_*^R(M)$ and
$\gamma\in F_*$, $x\.\gamma$ is defined as the composition:
\[
R \xra{\ x\ } F\wedge M \cong F\wedge R\wedge M \xra{F\wedge
\gamma \wedge M} F\wedge F\wedge M \xra{\mu\wedge M} F\wedge M
\]
The following fact is well-known for commutative $F$, see
\cite[Thm.\@ 4.2]{boardmanstable} for instance. The argument
generalizes to our situation. The point is that the smallest
localizing subcategory of $\dr$ \cite[Def.\@ 1.1.1]{hps} generated
by the suspensions of $R$ is $\dr$.

\begin{prop}\label{kunnethiso}
Assume that $F_*^R(M)$ or $F_*^R(N)$ is $F_*$--flat. Then the
pairing \eqref{kunneth} induces the {\em K\"unneth isomorphism} of
$F_*$--modules
\begin{equation*}\label{kunnethformula}
F_*^R(M)\otimes F_*^R(N) \cong F_*^R(M\wedge N).
\end{equation*}
\end{prop}

In the usual way, we construct a Kronecker pairing
\[
F_*^R(M)\otimes F^*_R(M) \lra F_*
\]
of left $F_*$--modules. We use the same right $F_*$-action on
$F_*^R(M)$ as in \eqref{kunneth} to form the tensor product.
Adjoint to the Kronecker pairing is a duality morphism
\begin{equation}\label{Rduality}
F^*_R(M) \lra \Hom_{F_*}(F_*^R(M), F_*)
\end{equation}
of left $F_*$-modules. We abbreviate $\Hom_{F_*}(-, F_*)$ by
$D(-)$ in the following.

In order to obtain K\"unneth isomorphisms in cohomology, we need
to consider a completed version of cohomology, analogous to the
one considered in \cite{boardmanstable}. First, we need an
analogue of the profinite filtration of cohomology groups as used
in \cite[Def.\@ 4.9]{boardmanstable}. Our objects, the
$R$--modules, are not naturally built from cells, but we can
approximate them by such objects. Recall from \cite[III.2]{ekmm}
that a cell $R$--module is a strict $R$--module $M$ which comes
with a certain filtration, called the sequential filtration. It
consists of sub $R$--mo\-dules $M_n$ with $M_0=\ast$ and with the
property that $M_{n+1}$ is obtained from $M_n$ by attaching a
wedge of $R$--spheres. Together with the $R$--linear maps which
respect this filtration, the cell $R$--modules form a subcategory
of $\mathscr{M}_R$. The approximation theorem \cite[Thm.\@
III.2.10]{ekmm} states that for each strict $R$--module $M$ there
exists a cell $R$-module $\Gamma M$ and a weak equivalence
$\gamma\: \Gamma M\to M$. The pair $(\Gamma M, \gamma)$ is called
a cell approximation of $M$. With this background, we are in a
position where we can mimick the definitions from
\cite{boardmanstable}.

\begin{defn}
Let $M$ be a cell $R$--module. The profinite filtration of
$F^*_R(M)$ consists of all the ideals
\[
F^a F^*_R(M) = \ker(F^*_R(M) \lra F^*_R(M_a))
\]
where $M_a$ runs through all finite cell submodules of $M$. These
ideals form a directed system. For an arbitrary $R$--module $M$,
choose a cell approximation $(\Gamma M, \gamma)$. The profinite
filtration on $F^*_R(\Gamma M)$ induces a filtration on $F^*_R(M)$
via the isomorphism
\[
F^*_R(\gamma)\: F^*_R(M) \cong F^*_R(\Gamma M).
\]
Define the {\em profinite topology\/} on $F^*_R(M)$ to be the
induced $F^*$--linear topology.
\end{defn}

As a consequence of \cite[Lemma III.2.2]{ekmm}, the profinite
topology is independent of the choice of a cell approximation.

\begin{defn}
For an $R$--module $M$, {\em completed $F$--cohomology\/}
$F^*_R(M)\hat\ $ is defined as the completion of $F^*_R(M)$ with
respect to the profinite topology.
\end{defn}

\begin{defn}
Let $M_*$ be an $F_*$--module. The dual-finite filtration on
$DM_*$ consists of the submodules
\[
F^L(DM_*) = \ker(DM_* \lra DL_*),
\]
where $L_*$ runs through all finitely generated submodules of
$M_*$. It gives rise to the {\em dual-finite topology\/} on
$DM_*$.
\end{defn}

The duality morphism \eqref{Rduality} is continuous with respect
to the profinite and the dual-finite topologies respectively. We
can prove the following statement in the same way as \cite[Thm.\@
4.14]{boardmanstable}.

\begin{prop}\label{dual}
Let $M$ be an $R$--module for which $F_*^R(M)$ is $F_*$--free.
Then the duality morphism induces a homeomorphism
\[
F^*_R(M) \cong D(F_*^R(M)).
\]
\end{prop}

In cohomology, we have an exterior product
\begin{equation}\label{kunnethcoh}
F^*_R(M)\otimes F^*_R(N) \lra F^*_R(M\wedge N),
\end{equation}
defined as usual. Again, we have to be careful with the right
action of $F_*$. For $x\in F^*_R(M)$ and $\gamma\in F_*$,
$x\cdot\gamma$ is defined as
\[
M\xra{\ x\ } F\cong F\wedge R \xra{F\wedge\gamma} F\wedge F\xra{\
\mu\ } F.
\]
As cell structures on $M$ and on $N$ induce one on $M\wedge N$
\cite[Prop.\@ III.7.3]{ekmm}, the product \eqref{kunnethcoh} lifts
to a pairing
\begin{equation}\label{kunnethcohcomplete}
F^*_R(M) \hatotimes F^*_R(N) \lra  F^*_R(M\wedge N)\hat\ .
\end{equation}
Here, the left hand side is the completion with respect to the
canonical topology on the tensor product. The following
K\"unneth-type theorem can be proved in the same way as
\cite[Thm.\@ 4.19]{boardmanstable}.

\begin{prop}\label{kunnethcohiso}
Assume that $F_*^R(M)$ is $F_*$--free. If $F_*^R(M)$ is finitely
generated or if $F_*^R(N)$ is $F_*$--free, there is an isomorphism
of left $F^*$-modules
\[
F^*_R(M\wedge N) \cong F^*_R(M) \hatotimes F^*_R(N).
\]
\end{prop}

\section{Regular triples of ring spectra}\label{regulartriples}

Let $M$ be an $R$--module. Recall the quotient constructions $M/x$
and $M/S$ from Section \ref{notation}, defined for an element $x$
and a sequence of elements $S$ in $R_*$ respectively.

\begin{defn}\label{defregquot}
If $S=(x_1, x_2, \ldots)$ is regular on $M_*$, we call $M/S$ a
{\em regular quotient\/} of $M$. Explicitly, this means that
multiplication by $x_i$ is injective on $M_*/(x_1, \ldots,
x_{i-1})M_*$ and that $M_*/SM_*\neq 0$.
\end{defn}

\begin{rem}
If $S$ and $S'$ are regular sequences generating the same ideal
$I$, the regular quotients $M/S$ and $M/S'$ are isomorphic
\cite[Cor.\@ V.2.10]{ekmm}. So it suffices to specify an ideal $I$
and set $M/I=M/S$ for some regular sequence $S$ generating $I$.
\end{rem}

Let $F$ be some $R$--ring spectrum and let $N=M/S$ be a quotient
of $M$ with respect to a sequence $S=(x_0, x_1, \ldots)$ in $R_*$.
Let $\Lambda_{F_*}(\alpha_0,\alpha_1, \ldots)$ denote the exterior
algebra on a set of generators $\alpha_j$ in bijective
correspondence with the $x_j$, with degrees $|\alpha_j|=|x_j|+1$.
Its $F_*$--dual is the completed exterior algebra
$\widehat\Lambda_{F_*}(\alpha_0^\vee, \alpha_1^\vee, \ldots)$ on
the dual basis $\{\alpha_j^\vee\}$ of $\{\alpha_j\}$. The
statement below follows from Propositions \ref{kunnethiso} and
\ref{kunnethcohiso}.

\begin{lem}\label{homregular}
Assume that $F_*^R(M)$ is $F_*$--free and that the elements $x_j$
are contained in the annihilator of $F_*$. Then there are
isomorphisms of $F_*$--modules:
\begin{align*}
F_*^R(N) & \cong F_*^R(M)\otimes F_*^R(\bigwedge_j R/x_j) \cong
F_*^R(M)\otimes \Lambda_{F_*}(\alpha_0, \alpha_1, \ldots)\\
F^*_R(N) & \cong F^*_R(M) \hatotimes F^*_R(\bigwedge_j R/x_j)
\cong  F^*_R(M) \hatotimes \widehat\Lambda_{F^*}(\alpha_0^\vee,
\alpha_1^\vee, \ldots)
\end{align*}
\end{lem}

Let now $(x_0, x_1, \ldots)$ be a regular sequence on $R_*$
consisting of non-zero divisors. Assume that $R$ is even, by which
we mean that the homotopy $R_*$ is concentrated in even degrees.
Then Strickland proves in \cite[Prop.\@ 3.1]{strickland} that the
$R$--modules $R/x_j$ can be realized as $R$--ring spectra. Fix
such products. Then there is a unique $R$--ring structure on $L =
\bigwedge_j R/x_j$ such that the maps $R/x_j \to L$ are commuting
ring maps, in the sense of \cite[Def.\@ 4.1]{strickland}.

\begin{defn}
Let $R$ be even. A {\em regular $R$--ring} is a regular quotient
of $R$ with respect to some regular sequence on $R_*$ consisting
of non-zero divisors, together with some chosen $R$--ring
structure arising as above.
\end{defn}

Let $L$ be a regular $R$--ring. Smashing a map from
$(R/x_j)^*_R(R/x_j)$ with the identities on the other smash
factors yields a homomorphism of $R_*$--algebras
\begin{align}\label{includej}
q_j\:  (R/x_j)^*_R(R/x_j) \lra L^*_R(L).
\end{align}
Lemma \ref{homregular} implies that $(R/x_j)^*_R(R/x_j)$ and
$L^*_R(L)$ are free over $R^*/(x_j)$ and $L^*$ respectively. Hence
they admit coalgebra structures. Clearly, the $q_j$ are coalgebra
maps. Let $\bar\beta_j$ be the map $\bar\beta_{x_j}$ as defined in
\eqref{standardbockstein}. Denote the image $q_j(\bar\beta_j) \in
L^*_R(L)$ by $Q_j$.

\begin{prop}\label{essential}
Let $R$ be even. If $L$ is a regular $R$--ring, there is an
isomorphism of $L_*$--modules
\[
L_*^R(L) \cong \Lambda_{L_*}(a_0, a_1, \ldots)
\]
and an isomorphism of $L_*$--bialgebras
\[
L^*_R(L) \cong \widehat\Lambda_{L^*}(Q_0, Q_1, \ldots)
\]
where $|a_j|=|Q_j|=|x_j|+1$.
\end{prop}

\begin{proof}
It follows from Lemma \ref{homregular} that the module structures
are as asserted. For the identification of $L^*_R(L)$ as an
algebra see \cite[Prop.\@ 4.15]{strickland}. The coalgebra
structure is determined by the fact that the $Q_j$ are derivations
\cite[Cor.\@ 4.19]{strickland} and hence are primitive with
respect to the coproduct.
\end{proof}

\begin{defn}\label{regulartriple}
Suppose that $R$ is an even commutative $\ss$--algebra and that
$T$ is a commutative $R$--ring spectrum. Assume that $F$ is an
$R$--ring spectrum of the form $F=T\wedge L$, where $L$ is a
regular $R$--ring with respect to a sequence $S=(x_0, x_1,
\ldots)$ in $R_*$ which is regular both on $R_*$ and on $T_*$.
Then we call $(R, T, F)$ a {\em regular triple}.
\end{defn}

We write $J=(x_0, x_1, \ldots)\lhd R_*$ and $I=J\cdot T_*\lhd T_*$
for the ideals generated by $S$ in the rings $R_*$ and $T_*$
respectively. Thus we have $L_*\cong R_*/J$ and $F_*\cong T_*/I$.
We abuse notation and write $x_j$ for the image in $T_*$ of the
elements $x_j\in R_*$ under $(\eta_T)_*\: R_*\to T_*$, induced by
the unit $\eta_T\: R\to T$ of $T$.

\begin{rem}\label{commutative}
If $2$ is a unit in $L_*$, $L$ may be realized as a commutative
$R$--ring spectrum \cite[Thm.\@ 4.11]{strickland}. So $F$ is then
commutative as an $R$--ring spectrum.
\end{rem}

The following is clear from the statements in Section
\ref{kunnethsec} and Proposition \ref{essential}. Because the maps
$L_*^R(L)\to F_*^R(L)$ and $L^*_R(L)\to F^*_R(L)$ induce
isomorphisms on tensoring the sources with $F_*$, we don't
introduce new labels for the generators of $F_*^R(L)$ and
$F^*_R(L)$. The condition that $F_*^R(T)$ is $F_*$--free is always
satisfied if $T$ is a regular quotient of $R$, as a consequence of
Lemma \ref{homregular}.

\begin{lem}
Let $(R, T, F)$ be a regular triple. Then there is an isomorphism
of $F_*$--modules
\[
F_*^R(F)\cong F_*^R(T)\otimes \Lambda_{F_*}(a_0, a_1, \ldots).
\]
If the sequence $S$ generating $I$ is finite or if $F_*^R(T)$ is
$F_*$--free, we have
\[
F^*_R(F) \cong F^*_R(T)\hatotimes \Lambda_{F^*}(Q_0, Q_1, \ldots)
\]
as $F^*$--modules.  In this case, the exterior algebra $F^*\otimes
L^*_R(L)$ maps isomorphically onto the subalgebra
$\Lambda_{F^*}(Q_0, Q_1, \ldots)$ under the natural map of
$F^*$--algebras
\[
F^*\otimes L^*_R(L) \lra F^*_R(F).
\]
\end{lem}

Recall the various $MU$--modules $v_n^{-1}BP$, $P(n)$, $B(n)$,
$BP\langle n\rangle$, $E(n)$, $\widehat{E(n)}$, $k(n)$, $K(n)$
derived from the Brown-Peterson spectrum $BP$ for some prime $p$,
with
\[
BP_* \cong \Z_{(p)}[v_1, v_2, \ldots]
\]
where $|v_n|=2(p^n-1)$ (see \cite{ravenelgreen} for instance). Let
$I_n=(v_0, \ldots, v_{n-1})$, where as usual $v_0=p$. We have
isomorphisms of $BP_*$--modules
\begin{align*}
P(n)_* & \cong \F_p[v_n, v_{n+1}, \ldots] & B(n)_* & \cong
v_n^{-1} P(n)_*
\\
BP\langle n\rangle_* & \cong \Z_{(p)}[v_1, \ldots, v_n] & E(n)_* &
\cong v_n^{-1} BP\langle n\rangle_*
\\
(v_n^{-1}BP)_* & \cong v_n^{-1} BP_* & \wh{E(n)}_* & \cong
(E(n)_*)^\wedge_{I_n}
\\
k(n)_* & \cong \F_p[v_n] & K(n)_* & \cong v_n^{-1} k(n)_*
\end{align*}
where $(-)^\wedge_{I_n}$ is completion at $I_n$. Using results
from \cite{strickland}, we find the following examples of regular
triples:

\begin{exmp}\label{classicalmu}
$R=MU$, $T=E(n)$ and $F=K(n)$ for an odd prime $p$. For $p=2$, the
substitute $E(n)'$ is commutative, see \cite{strickland}. Hence
$(MU, E(n)', K(n))$ is regular.
\end{exmp}

\begin{exmp}
From \cite{br} or \cite{goersseinfty} we know that $\wh{E(n)}$ is
a commutative $\ss$--algebra. So $(\wh{E(n)}, \wh{E(n)}, K(n))$ is
a regular triple. Of course, $(MU, \wh{E(n)}, K(n))$ is regular as
well.
\end{exmp}

\begin{exmp}\label{examplebp}
$R=MU$, $T=BP$, $F=P(n)$.
\end{exmp}

\begin{exmp}
$R=MU$, $T=P(n)$, $F=k(n)$ for odd $p$.
\end{exmp}

\begin{exmp}
$R=MU$, $T=BP\langle n\rangle$, $F=k(n)$ for odd $p$. For $p=2$,
we may take $BP\langle n\rangle'$ instead, see \cite{strickland}.
\end{exmp}

\begin{exmp}
$R=v_n^{-1} MU$, $T=v_n^{-1} BP$, $F=K(n)$.
\end{exmp}

\begin{exmp}
$R=v_n^{-1} MU$, $T=B(n)$, $F=K(n)$ for odd $p$.
\end{exmp}

We end the section by deriving a characterization of the morphisms
$L_*^R(Q_j)$ induced by the Bockstein operations on $L$--homology
groups. By definition, the morphism induced by $Q_j$ on
$L$--cohomology
\[
L^*_R(Q_j)\:  L^*_R(L) \lra L^{*+|x_j|+1}_R(L)
\]
is multiplication by $Q_j$ on the right. In the following $D(-)$
denotes the continuous $L_*$--duality operator. Consider
\begin{equation}\label{qcohdual}
D(L^*_R(Q_j)) \: D(L^{*+|x_j|+1}_R(L)) \lra D(L^*_R(L)).
\end{equation}
The coalgebra structure on $L^*_R(L)$ induces an algebra structure
on $D(L^*_R(L))$, which is just the completed exterior algebra
generated by the dual basis $\{Q_j^\vee\}$ of $\{Q_j\}$. Because
$L^*_R(Q_j)$ is multiplication by a primitive element, the dual
map $D(L^*_R(Q_j))$ is a derivation with respect to this algebra
structure. As such, it is determined by
\[
D(L^*_R(Q_j))(Q_i^\vee) = \delta_{ij},
\]
where $\delta_{ij}$ is the Kronecker delta. Now as $L_*^R(L)$ is a
free $L_*$--module, $L_*^R(L)$ injects into its bidual
$DD(L_*^R(L))$ and $DD(L_*^R(Q_j))$ restricts to $L_*^R(Q_j)$. Via
the duality isomorphism
\[
L^*_R(L) \cong D(L_*^R(L)),
\]
the morphism $D(L_*^R(Q_j))$ corresponds to $L^*_R(Q_j)$, and
therefore \eqref{qcohdual} describes $L_*^R(Q_j)$. Hence, choosing
the $a_j\in L_*^R(L)$ and their formal exterior products so that
they correspond to the $Q_j^\vee$ and their products, we have
shown the following.

\begin{lem}\label{bocksteinhomology}
Endow $L_*^R(L) \cong \Lambda_{L_*}(a_0, a_1, \ldots)$ with the
exterior algebra structure. Then the homomorphism
\[
L_*^R(Q_j) \: L_*^R(L) \lra L_{*-|x_j|-1}^R(L)
\]
is the derivation $\frac{\partial}{\partial a_j}$ characterized by
$\frac{\partial}{\partial a_j}(a_i) = \delta_{ij}$.
\end{lem}

\begin{rem}\label{commtwoalgebra}
If $L$ is commutative, the algebra structure of $L_*^R(L)$ induced
by the multiplication on $L$ is the one considered in the lemma.
This follows from the fact that the duality map $d\: L^*_R(L)
\cong D(L_*^R(L))$ is an algebra map under this assumption. For
non-commutative $L$, the two algebra structures will not agree in
general.
\end{rem}

\section{The algebra of $I$--adic towers}\label{algebraiadic}

Let $(R, T, F)$ be a regular triple with $F\cong T\wedge L$.
Recall the isomorphisms $L_* \cong R_*/J$ and $F_* \cong T_*/I$,
where $J$ is generated by a regular sequence $(x_0, x_1, \ldots)$
and $I=J\cdot T_*$. The contents of this section apply to both
$(R_*, J_*)$ and $(T_*, I_*)$. We formulate everything for the
first case and write $\otimes$ for $\otimes_{L_*}$ in this
section.

Define two families $(\cale^s)_{s\geq 0}$ and $(\calf^s)_{s\geq
0}$ of short exact sequences of $R_*$-modules
\begin{equation*}\xymatrix@R=0.3cm{ \cale^s : &   0 \ar[r] &
J^{s+1} \ar[r]^{i_{s+1}} & J^s \ar[r]^-{p_{s}} & J^s/J^{s+1}
\ar[r]& 0\\ \calf^s : &   0 \ar[r] & J^s/J^{s+1} \ar[r]^{j_{s}}
& R_*/ J^{s+1} \ar[r]^-{q_{s+1}} & R_*/J^s  \ar[r] & 0}
\end{equation*}
where by convention $J^0=R_*$. All the maps are the canonical
ones. Whenever it is clear from the context, we will omit the
index $s$. We will also write $i$ and $q$ for a composition of
maps $i$ and $q$ as defined above.

Let $\mathscr D_{R_*}$ be the derived category of the ring $R_*$.
Its objects are the chain complexes of graded $R_*$--modules. The
sequences $\cale^s$ and $\calf^s$ induce towers in $\mathscr
D_{R_*}$ consisting of triangles, the {\em $J$--adic tower over
$R_*$}
\begin{equation}
\label{algebraiciadicover} \begin{array}{c} \xymatrix@!C=.8cm{
\cdots\ \ar[r]^{i} & J^3 \ar[r]^{i}\ar[d]^{p} & J^2
\ar[r]^{i}\ar[d]^{p} & J \ar[r]^{i}\ar[d]^{p} & R_*\ar[d]^{p} \\
& J^3/J^4\ar[lu]^{\epsilon^3}|-\circ &
J^2/J^3\ar[lu]^{\epsilon^2}|-\circ &
J/J^2\ar[lu]^{\epsilon^1}|-\circ & L_*\ar[lu]^{\epsilon^0}|-\circ
}
\end{array}
\end{equation}
and the {\em $J$-adic tower under $R_*$}:
\begin{equation}
\label{algebraiciadicunder} \begin{array}{c} \xymatrix@!C=1cm{ R_*
\ar[r] & \ \cdots\ \ar[r]^-{q} & R_*/J^3 \ar[r]^-{q}
\ar[ld]^{\theta^3}|-\circ & R_*/J^2 \ar[r]^-{q}
\ar[ld]^{\theta^2}|-\circ & L_* \ar[ld]^{\theta^1}|-\circ
\\
& \cdots & J^2/J^3 \ar[u]_{j} & J/J^2 \ar[u]_{j}}
\end{array}
\end{equation}
Arrows with a circle denote morphisms of degree $-1$.

Let $M_*$ be a graded $R_*$--module. We write $\gr^*_J(M_*)$ for
the (bi)graded $L_*$--module associated to the $J$--adic
filtration on $M_*$. For $M_*=R_*$, we have
\[
\gr^*_J(R_*) = \bigoplus_{s\geq 0} J^s/J^{s+1}.
\]
For a graded $L_*$--module $Y_*$, $\Sym^*_{L_*}(Y_*)$ denotes the
symmetric algebra over $L_*$ on $Y_*$. It is naturally bigraded.

\begin{prop}[{\cite[Thm.\@ 16.2]{matsumura}}]
\label{splitting} The $L_*$--module $J/J^2$ is freely generated by
the residue classes $\{x_j\}$ of the generators $x_j\in J$.
Furthermore, we have an isomorphism of graded $L_*$--algebras
\[
\gr^*_J(R_*) \cong \Sym^*_{L_*}(J/J^2).
\]
\end{prop}

It follows that we have isomorphisms of $L^*$--modules
\begin{equation}\label{extgr}
\Ext^{*,*}_{R^*}(L^*, J^s/J^{s+1}) \cong \Ext^{*,*}_{R^*}(L^*,
L^*)\otimes J^s/J^{s+1}.
\end{equation}
The element of $\Ext^{1,0}_{R^*}(L^*, J/J^2)$ associated to the
short exact sequence
\[
\calf^1: \quad 0 \lra J/J^2 \lra R_*/J^2 \lra L_* \lra 0
\]
corresponds via \eqref{extgr} to a sum of the form $\sum
f_j\otimes \{x_j\}$ with
\[
f_j\in\Ext^{1, -|x_j|}_{R^*}(L^*, L^*).
\]
The following fact is as well-known as fundamental, see
\cite[Thm.\@ 16.5]{matsumura} and for the multiplicative structure
\cite[\S 4.5]{weibel}. The identification of the elements $f_j$ as
generators is straight-forward.

\begin{prop}\label{exterior}
There are isomorphisms of $L_*$--algebras
\begin{align*}
\Ext^{*,*}_{R^*}(L^*, L^*)  & \cong \widehat\Lambda_{L^*}(f_0,
f_1, \ldots)\\
\Tor_{*,*}^{R_*}(L_*, L_*) & \cong \Lambda_{L_*}(e_0, e_1,
\ldots)
\end{align*}
where the $f_j$ and the $e_j$ have bidegree $(1, -|x_j|)$ and $(1,
|x_j|)$ respectively.
\end{prop}

Let $\derotimes$ and $\RHom$ be the left and right derived
functors of $\otimes_{R_*}$ and $\Hom_{R_*}$ respectively. The
functors
\[
HL_{p,*}(-) = H_p(L_*\derotimes-), \quad HL^{p,*}(-) =
H^p(\RHom(-, L^*))
\]
define a homology and a cohomology theory on $\mathscr{D}_{R_*}$
respectively. We call them $L_*$--homology and $L_*$--cohomology.
Both $HL_{*,*}(-)$ and $HL^{*,*}(-)$ take values in the category
of $L_*$--modules. Their restriction to $R_*$--modules are the
functors $\Tor_{*,*}^{R_*}(L_*, -)$ and $\Ext^{*,*}_{R^*}(-, L^*)$
respectively.

It is well-known that the cohomology theory $E^*(-)$ takes values
in $E^*(E)$--modules if $E$ is a ring spectrum. Also, if $E_*(E)$
is $E_*$--flat, $E_*(E)$ has a natural coalgebra structure and
$E_*(-)$ takes values in $E_*(E)$--comodules. Mimicking the
definition of these structures, we can define a natural coproduct
on $HL_{*,*}(L_*, L_*)$ and show that the functors $HL_{*,*}(-)$
and $HL^{*,*}(-)$ take values in the categories of
$HL_{*,*}(L_*)$--comodules and $HL^{*,*}(L^*)$--modules
respectively. Arguing as in Remark \ref{commtwoalgebra}, we find
that the product on $HL^{*,*}(L^*)$ is dual to the coproduct on
$HL_{*,*}(L_*)$. Hence the latter agrees with the natural
coproduct on $\Lambda_{L_*}(e_0, e_1, \ldots)$.

The following result, proved as \cite[Thm.\@ 3.8]{swglasgow} or
\cite[Thm.\@ 1]{swshort}, describes how surprisingly simple the
$L_*$--homology of the $J$--adic tower over $R_*$ is. In
\cite{swglasgow} and \cite{swshort}, we only considered finite
regular quotients. The general case, however, follows easily by
passing to limits.

\begin{thm}\label{homtower}
The functor $HL_{*,*}(-)$ maps the algebraic $J$--adic tower over
$R_*$ to a diagram of $HL_{*,*}(L_*)$--comodules composed of short
exact sequences of the form
\begin{equation} \label{khomtower}
\begin{array}{c}
\xymatrix{L_*  \ar[d]^{p_*} & HL_{*,*}(J)\ar[l]_-{0} \ar[d]^{p_*}&
HL_{*,*}(J^2)\ar[l]_-{0}\ar[d]^{p_*} &\ \cdots. \ar[l]_-0
\\
HL_{*,*}(L_*)\ar[ur]_{\epsilon^0_*}|-\circ &
HL_{*,*}(J/J^2)\ar[ur]_{\epsilon^1_*}|-\circ & HL_{*,*}(J^2/J^3)
\ar[ur]_{\epsilon^2_*}|-\circ & \ \cdots}
\end{array}
\end{equation}
Furthermore, all modules are free over $L_*$.
\end{thm}

From this, we easily derive the $L_*$--homology of the $J$--adic
tower under $R_*$. Namely, we have an isomorphism
$HL_{*,*}(R_*/J^s) \cong L_* \oplus HL_{*-1,*}(J^s)$, induced by
the long exact sequence of homology groups associated to
\begin{equation}\label{sesses}
0 \lra J^s \xrightarrow{\ i\ } R_* \xrightarrow{\ q\ } R_*/J^s
\lra 0.
\end{equation}
$$
\wt{HL}_{*,*}(R_*/J^s) = \ker(HL_{*,*}(R_*/J^s) \to L_*) \cong
HL_{*-1,*}(J^s).\leqno{\hbox{We set}}
$$

\begin{cor}\label{homquottower}
Applying $HL_{*,*}(-)$  to the algebraic $J$--adic tower under
$R_*$ yields an unrolled exact couple of the form:
\[
\xymatrix{\wt{HL}_{*,*}(L_*) \ar[d]^{\theta^1_*}|-\circ&
\wt{HL}_{*,*}(R/J^2) \ar[l]_-0\ar[d]^{\theta^2_*}|-\circ &
\wt{HL}_{*,*}(R/J^3) \ar[l]_-0\ar[d]^{\theta^3_*}|-\circ &\
\cdots \ar[l]_-0
\\
HL_{*,*}(J/J^2) \ar[ur]_-{j_*}& HL_{*,*}(J^2/J^3) \ar[ur]_-{j_*} &
HL_{*,*}(J^3/J^4)\ar[ur]_-{j_*} & \ \cdots}
\]
Furthermore, all modules are free over $L_*$.
\end{cor}

\begin{proof}
Let $\phi^s$ be the connecting homomorphism associated to
\eqref{sesses}. By naturality of connecting homomorphisms, we have
a commutative diagram:
\[
\xymatrix{ HL_{*-1,*}(J^{s+1})\ar[r]^-{i_*} & HL_{*-1,*}(J^s)
\\
\wt{HL}_{*,*}(R/J^{s+1})\ar[r]^-{q_*} \ar[u]^{\cong}_-{\phi^{s+1}}
& \wt{HL}_{*,*}(R/J^s) \ar[u]^{\cong}_-{\phi^s} }
\]
Hence the theorem implies that $\wt{HL}_{*,*}(R/J^{s+1})
\xrightarrow{q_*} \wt{HL}_{*,*}(R/J^s)$ is trivial.
\end{proof}

\begin{rem}\label{cohomologydual}
As all homology groups appearing in Theorem \ref{homtower} and
Corollary \ref{homquottower} are $L_*$--free, the
$L_*$--cohomology of the algebraic $J$--adic towers is isomorphic
to the $L_*$--dual of their respective $L_*$--homology. This
follows from a statement analogous to Proposition \ref{dual}, but
for the algebraic derived category $\mathscr D_{R_*}$.
\end{rem}

By Proposition \ref{splitting}, there are isomorphisms of
$L_*$--modules
\begin{equation}\label{cofree}
HL_{*,*}(J^s/J^{s+1}) \cong HL_{*,*}(L_*) \otimes J^s/J^{s+1}.
\end{equation}
In fact, these are isomorphisms of $HL_{*,*}(L_*)$--comodules, if
we regard the right hand side as the cofree
$HL_{*,*}(L_*)$--comodule generated by $J^s/J^{s+1}$. As all
modules in \eqref{khomtower} are $L_*$--free, we obtain a long
exact sequence with the following property (see \cite[Ch.\@
IX]{maclanehomology} for a discussion of relative homological
algebra):

\begin{cor}\label{relinjalg}
The sequence of $HL_{*,*}(L_*)$--comodules
\[
0 \to L_* \xrightarrow{p_*} HL_{*,*}(L_*) \xrightarrow{\delta^0}
HL_{*-1,*}(L_*)\otimes J/J^2 \xrightarrow{\delta^1}
HL_{*-2,*}(L_*)\otimes J^2/J^3 \to \cdots,
\]
where $\delta^s = p_*\epsilon^s_*$, is a relative injective
resolution of $L_*$.
\end{cor}

We may consider the differentials $\delta^s$ of this complex as
the components of a single differential
\begin{equation}\label{gradeddiff}
\delta^*\: HL_{*,*}(L_*)\otimes \gr_J^*(R_*) \lra
HL_{*-1,*}(L_*)\otimes \gr_J^{*+1}(R_*)
\end{equation}
on the $\gr_J^*(R_*)$--algebra $HL_{*,*}(L_*)\otimes
\gr_J^*(R_*)$. By \cite[Prop.\@ 3.10]{swglasgow} or \cite[Prop.\@
7]{swshort}, it can be characterized as follows.

\begin{prop}\label{connectinghom}
The differential $\delta^*$ of \eqref{gradeddiff} is a derivation
over $\gr_J^*(R_*)$. As such, it is determined by its values
$\delta^*(e_j)= - \{x_j\}$ on the generators $e_j$ of the algebra
$HL_{*,*}(L_*)$.
\end{prop}

Explicitly, for an element $e_{i_1}\wedge\dotsb\wedge
e_{i_l}\otimes \bar f \in HL_{*,*}(L_*)\otimes \gr_J^*(R_*)$, we
have
\[
\delta^*(e_{i_1}\wedge\dotsb\wedge e_{i_l}\otimes \bar f) =
\sum_{j=1}^{l} (-1)^{j+1} \ e_{i_1}\wedge\dotsb\wedge \widehat
e_{i_j}\wedge\dotsb \wedge e_{i_l}\otimes \overline{x_{i_j} f},
\]
where the hat indicates that the entry underneath is omitted.

\section{Adams resolutions}\label{adamsresolutions}

In this section, we recall the notion of Adams resolutions with
respect to a given injective class in a triangulated category. We
then indicate how the injective classes giving rise to the
traditional Adams resolutions can be obtained in a natural way. In
an analogous way, we obtain the injective class we need for the
construction of the $I$--adic tower. We are rather brief here. A
more detailed exposition can be found in \cite{chr}, see also
\cite{eilenbergmoore} and \cite{husemollermoore}. The original
reference for Adams resolutions (in the classical sense) is
\cite{miller}.

A {\em length two complex\/} in an additive category $\calc$ is a
diagram
\begin{equation}\label{allex}
F \xra{\ f\ } L \xra{\ g\ } M
\end{equation}
with $gf=0$. A {\em complex\/} is a sequence
\[
\cdots \lra M_{s+1} \xra{\ f_{s+1}\ } M_s \xra{\ f_s\ } M_{s+1}
\lra \cdots
\]
such that $f_s f_{s+1}=0$ for all $s$. If the terms $M_s$ of a
complex are zero for all large $s$ or all small $s$, we omit them
from the notation. An {\em injective class\/} $\calj$ in $\calc$
is specified by a collection of length two complexes, the
so-called {\em allowable exact complexes}, and a collection of
objects called the {\em injectives}. The two collections are
required to determine each other in the following way. An object
$J$ is injective if and only if the sequence
\[
[M, J]^* \xrightarrow{\ g^*\ } [L, J]^* \xrightarrow{\ f^*\ } [F,
J]^*
\]
of graded abelian groups is exact for each allowable exact complex
of the form \eqref{allex}, and vice versa. Furthermore, one
requires that there are enough injectives, in the following sense.
For each map $f\: M \to N$ there exists an injective $I$ and a map
$g\: N \to I$ such that $M\xra{f} N\xra{g} I$ is an allowable
exact complex.

Now assume that $\calc$ is triangulated. Then we may specify
instead of the allowable exact complexes a family of morphisms
called the {\em null-maps}. An object $J$ is then injective if and
only if the homomorphism of graded abelian groups
\[
[N, J]^* \xrightarrow{f^*} [M, J]^*
\]
is zero for each null map $f: M \to N$, and vice versa.

An example for an injective class is the {\em split injective
class}, defined in each additive category $\calc$. Let $\vee$
denote the biproduct (coproduct and product) in $\calc$. The
allowable complexes are the ones of the form
\[
A \vee B \lra B\vee C \lra C \vee D
\]
where the maps are the projections onto the second summand
followed by the inclusions into the first summand. Every object is
injective with respect to this injective class. If $\calc$ is
triangulated, the null-maps are the trivial maps.

There is the dual concept of a {\em projective class\/} in an
additive category $\calc$. It is specified by a collection of
length two complexes called the {\em allowable exact complexes},
and a collection of objects called the {\em projectives}. An
object $P$ is projective if and only if
\[
[P, F]^* \xrightarrow{f_*} [P, L]^* \xrightarrow{g_*} [P, M]^*
\]
is exact for each allowable exact complex of the form
\eqref{allex}, and vice versa. There have to exist enough
projectives, in an obvious sense. If $\calc$ is triangulated, we
can specify instead of the allowable exact sequences the
collection of null-maps. They are all the maps $f\: M\to N$ for
which
\[
f_*\: [P, M]^* \lra [P, N]^*
\]
is trivial for each projective $P$, and vice versa.

Suppose that $\calj$ is an injective class in an additive category
$\calc$. A {\em relative injective resolution\/} of an object $M$
with respect to $\calj$ is a complex of the form
\begin{equation}\label{allowableresolution}
0 \lra M \lra J_0 \lra J_1 \lra J_2 \lra \cdots
\end{equation}
with injectives $J_i$ which has the property that each three-term
subcomplex
\[
J_{s-1} \lra J_s \lra J_{s+1},
\]
where $J_{-1}=M$ and $J_s=0$ for $s<-1$, is allowable exact. If
$\calc$ is triangulated, such a resolution gives rise to an {\em
Adams resolution\/} of $M$ with respect to $\calj$, which is
unique up to isomorphism. By this, we mean a diagram of the form
\begin{equation*}
\begin{array}{c}
\xymatrix@!=0.3cm{  M=M_0 \ar[dr]  && M_1 \ar[ll]\ar[dr]
&& M_2 \ar[ll]\ar[dr] && M_3 \ar[ll]&\ar@{}[d]|{.\,.\,.}\\
& I_0 \ar[ur] \ar[ur]|-\circ && I_1 \ar[ur]|-\circ \ar[ur] && I_2
\ar[ur]|-\circ & & }
\end{array}
\end{equation*}
consisting of triangles
\[
M_{i+1} \lra M_i \lra I_i \lra \Sigma M_{i+1}
\]
such that all $I_i$ are injectives and all maps $M_{i+1} \to M_i$
are null-maps. Here $\Sigma$ denotes the suspension of $\calc$.
The injectives that we get from a relative injective resolution
\eqref{allowableresolution} are $I_i = \Sigma^{-i} J_i$.
Conversely, we can derive from an Adams resolution a relative
injective resolution. Applying $[N, -]^*$ to an Adams resolution,
where $N$ is another object of $\calc$, yields an exact couple of
abelian groups. It gives rise to an {\em Adams spectral sequence}.

Dually defined are relative projective resolutions and Adams
resolutions with respect to a projective class.

There is a standard method for transferring an injective class
from one additive category to another, namely by means of an
adjunction. Let $\calc$ and $\calc'$ be additive categories.
Assume that we have an adjoint pair of additive functors
\[
\xymatrix{\calc\ \ar@<0.5ex>[r]^F & \ \calc', \ar@<0.5ex>[l]^-G}
\]
where $F$ is the left adjoint. Given an injective class $\calj'$
in $\calc'$, we define an injective class $F^{-1}(\calj')$ in
$\calc$ as follows. The allowable exact complexes are all length
two complexes whose image under $F$ is allowable exact in
$\calc'$. The associated injectives turn out to be the retracts of
objects of the form $G(X)$, where $X$ is a $\calj'$--injective.

The injective classes in the stable homotopy category $\ds$ which
give rise to classical Adams resolutions arise in this way.
Namely, for a ring spectrum $R$, we have an adjoint pair of
additive functors
\[
\xymatrix{\ds\ \ar@<0.5ex>[rr]^-{V=R\wedge_\ss -\ } && \ \ds.
\ar@<0.5ex>[ll]^-{U=F_\ss(R, -)}}
\]
The {\em Adams injective class with respect to $R$\/} is
$V^{-1}(\calj)$, where $\calj$ is the split injective class in
$\ds$. Its injectives are the retracts of spectra of the form
$F_\ss(R,X)$. These agree with the retracts of spectra of the form
$R\wedge_\ss X$ (see Remark \ref{adamsmonoid}). More generally, if
$F$ is a ring spectrum over a commutative $\ss$--algebra $R$, we
have an additive adjunction
\begin{equation}\label{adjrtf}
\xymatrix{ \dr\ \ar@<0.5ex>[rr]^-{V=F\wedge -} && \ \dr.
\ar@<0.5ex>[ll]^-{U=F_R(F, -)}}
\end{equation}

\begin{defn}\label{defadams}
Let $\calj$ be the split injective class in $\dr$. The {\em Adams
injective class associated to $F$ in $\dr$\/} is the injective
class given via \eqref{adjrtf} as $V^{-1}(\calj)$.
\end{defn}

Similarly to the situation in $\ds$, the injectives are the
retracts of $R$--modules of the form $F\wedge M$, where $M$ is an
$R$--module.

\begin{rem}\label{adamsmonoid}
Note that we have an adjunction \eqref{adjrtf} for any object $F$
in $\dr$, not just for ring spectra. It is the recognition of the
injectives mentioned above which depends on having a
multiplicative structure on $F$. Indeed, this implies that $V$ is
a monad and $U$ a comonad, and as $V$ and $U$ are adjoint, the
$V$--algebras (the retracts of $F\wedge M$, $M\in\dr$) are the
same as the $U$--coalgebras (the retracts of $F_R(F, M)$,
$M\in\dr$).
\end{rem}

\begin{rem}\label{adamsalgebraic}
In a completely analogous way, one can define an Adams injective
class associated to a monoid $F_*$ in the derived category
$\mathscr D_{R_*}$ of a graded ring $R_*$. In particular, this
applies to an algebra $F_*$ over $R_*$ or more generally to a
differential graded algebra $F_*$. Instead of \eqref{adjrtf}, one
uses the adjunction
\begin{equation}
\xymatrix{ \mathscr D_{R_*}\ \ar@<0.5ex>[rr]^-{V=F_*\derotimes -}
&& \ \mathscr D_{R_*} \ar@<0.5ex>[ll]^-{U=\RHom(F_*, -)}}
\end{equation}
(see Section \ref{algebraiadic} for the definition of $\derotimes$
and $\RHom$).
\end{rem}

We also mention the {\em ghost projective class}. It is defined in
the triangulated categories $\mathscr{D}_R$, where $R$ is an
$\ss$--algebra. The null-maps are the morphisms which are trivial
on homotopy groups. The projectives are the retracts of wedges of
suspensions of $R$.

\section{Construction of the tower}\label{topologytower}

Let $(R, T, F)$ be a regular triple with $F\cong T\wedge L$. Thus
we have $L_*\cong R_*/J$ and $F_*\cong T_*/I$, where $J$ is
generated by a sequence $S=(x_0, x_1, \ldots) \subset R_*$,
regular on $R_*$ and on $T_*$, and where $I=J\cdot T_*$. Consider
the graded polynomial ring $L_*[y_0, y_1, \ldots]$ with
$|y_j|=|x_j|$. Define $L$--free $R$--modules $J^s/J^{s+1}$ by
setting
\[
J^s/J^{s+1} = \bigvee_{y\in V_s} \Sigma^{|y|} L,
\]
where $V_s$ is the set of monomials of homogeneous degree $s$ in
the variables $y_j$. By ``$L$--free'', we mean isomorphic to a
wedge of suspensions of copies of $L$.

We now construct a sequence in $\dr$ of the form:
\begin{equation*}
\ast \lra R \xrightarrow{\ \eta\ } L \xrightarrow{\ \wt\epsilon^0\
} \Sigma J/J^2 \xrightarrow{\ \Sigma\wt\epsilon^1\ } \Sigma^2
J^2/J^3 \xrightarrow{\ \Sigma^2 \wt\epsilon^2\ } \cdots
\end{equation*}
The maps $\wt\epsilon^i$ are easy to describe. Namely, as all
terms $J^s/J^{s+1}$ are $L$--free, a map between them corresponds
to a matrix over the graded endomorphism ring $L^*_R(L)$. The
Bockstein operations define a map
\[
\wt\epsilon^0 = \bigvee_j Q_j \: L \lra \Sigma J/J^2.
\]
For $s>0$, we define $\wt\epsilon^s \: J^s/J^{s+1} \lra \Sigma
J^{s+1}/J^{s+2}$ to be the sum over all monomials $y\in V_s$ of
the maps
\[
\bigvee_j \Sigma^{|y|} Q_j \: \Sigma^{|y|} L \lra \bigvee_j
\Sigma^{|y y_j|+1} L.
\]
Recall from Proposition \ref{essential} that $L_*^R(L) \cong
\Lambda_{L_*}(a_0, a_1, \ldots)$ as an $L_*$--module.

\begin{lem}\label{relinjstatement}
The $L$--homology of the sequence
\begin{equation}\label{resolutionsequence}
\ast \lra R \xrightarrow{\ \eta\ } L \xrightarrow{\ \wt\epsilon^0\
} \Sigma J/J^2 \xrightarrow{\ \Sigma\wt\epsilon^1\ } \Sigma^2
J^2/J^3 \xrightarrow{\ \Sigma^2 \wt\epsilon^2\ } \cdots
\end{equation}
is a long exact sequence of $\Lambda_{L_*}(a_0, a_1,
\ldots)$--comodules
\begin{equation}\label{relinj}
0 \lra L_* \xrightarrow{\ \eta_*\ } L_*^R(L) \xrightarrow{\
\wt\epsilon^0_*\ } L_*^R(\Sigma J/J^2) \xrightarrow{\
\wt\epsilon^1_*\ } L_*^R(\Sigma J^2/J^3) \lra \cdots
\end{equation}
which provides a relative injective resolution of $L_*$.
\end{lem}

\begin{proof}
The isomorphism $L_*^R(L) \cong \Lambda_{L_*}(a_0, a_1, \ldots)$
from Proposition \ref{essential} induces
\[
L_*^R(J^s/J^{s+1}) \cong L_*^R(L)\otimes_{L_*} J^s/J^{s+1}.
\]
We have identified in Lemma \ref{bocksteinhomology} the morphism
$L_*^R(Q_j)$. Recall that we have obtained in Corollary
\ref{relinjalg} and Proposition \ref{connectinghom} a relative
injective resolution of $L_*$ over $\Lambda_{L_*}(e_0, e_1,
\ldots)$. It follows from the definition of the maps
$\wt\epsilon^s$ that mapping the generators $a_j$ of
$\Lambda_{L_*}(a_0, a_1, \ldots)$ to the generators $e_j$ yields
an isomorphism between \eqref{relinj} and the sequence in
Corollary \ref{relinjalg}, if we endow the latter with the total
gradation.
\end{proof}

From now on, we will omit the index $s$ and just write
$\wt\epsilon$ for all the maps in the sequence. To obtain an
abstract characterization of the sequence
\eqref{resolutionsequence} by means of Lemma
\ref{relinjstatement}, we need the following technical lemma.
Recall from Section \ref{notation} the difference between the
concepts of a module spectrum over a ring spectrum and a module
over an $\ss$--algebra. We write $\Cohom^*$ for the maps of graded
comodules over some graded coalgebra.

\begin{lem}\label{kmodulemaps}
Let $E$ be an arbitrary $R$--ring spectrum. Let $Y$ be an
$R$--module for which $E_*^R(Y)$ is $E_*$--free and let $Z$ be an
$E$--module spectrum in $\dr$. Then there is a natural isomorphism
of $E_*$--modules
\begin{equation*}\label{kduality}
[Y, Z]_R^* \cong \Hom_{E_*}^*(E_*^R(Y), Z_*).
\end{equation*}
If $E$ is a commutative $R$--ring spectrum and $E_*^R(E)$ is
$E_*$--flat, the functor $E_*^R(-)$ induces an isomorphism of
$E_*$--modules
\[
[Y, Z]_R^* \cong \Cohom^*_{E_*^R(E)}(E_*^R(Y), E_*^R(Z)).
\]
\end{lem}

\begin{proof}
First, we construct a natural transformation of functors of
$E$--module spectra with values in $E_*$--modules
\begin{equation}\label{first}
[Y, -]_R^* \lra \Hom^*_{E_*}(E_*^R(Y), (-)_*)
\end{equation}
as follows. The $E$--homology functor gives rise to a natural
transformation
\begin{equation}
[Y, -]_R^* \lra \Hom_{E_*}^*(E_*^R(Y), E_*^R(-)).
\end{equation}
Composing this with the transformation $E_*^R(-)\to (-)_*$ induced
by the action map for $E$--module spectra induces a transformation
as in \eqref{first}.

If we precompose both sides of \eqref{first} with $E\wedge_R-$, we
obtain a transformation
\begin{equation}\label{second}
[Y, E\wedge -]_R^* \lra \Hom_{E_*}^*(E_*^R(Y), E_*^R(-))
\end{equation}
of functors of $R$--modules. As $E_*^R(Y)$ is $E_*$--free, the
functor on the right hand side is exact. Hence both sides of
\eqref{second} define homology theories on $\dr$. On suspensions
of $R$, \eqref{second} is the duality morphism \eqref{Rduality},
and is therefore an isomorphism. Hence it is an equivalence on
$\dr$ by the usual argument. So \eqref{first} is an equivalence on
all $E$--module spectra $Z$ of the form $E\wedge X$, $X$ an
$R$--module. But an arbitrary $E$--module spectrum $Z$ in $\dr$ is
a retract of $E\wedge Z$, via the action map. So we obtain an
isomorphism for general $Z$ by naturality.

For commutative $E$ for which $E_*^R(E)$ is $E_*$--flat,
$E_*^R(-)$ takes values in $E_*^R(E)$--co\-mo\-dules. The
structure maps are defined exactly as in the case $R=S$. Consider
the composition of natural transformations of functors of
$E$--module spectra
\begin{equation}\label{third}
E_*^R(E) \otimes (-)_* \xra{E_*^R(E)\otimes\eta_*} E_*^R(E)\otimes
E_*^R(-) \xra{\ \times\ } E_*^R(E\wedge -) \lra E_*^R(-),
\end{equation}
where the last morphism is induced by the $E$--action. It can be
checked to be $E_*^R(E)$--colinear. Precomposing it with
$E\wedge-$ yields a transformation
\[
E_*^R(E) \otimes E_*^R(-) \lra E_*^R(E\wedge -)
\]
of homology theories on $\dr$. It does {\em not\/} agree with the
one induced by the exterior product $\times$. Nevertheless, an
inspection shows that it defines an equivalence of homology
theories. With an argument as above, we conclude that
\eqref{third} is an equivalence. Now recall that we have an
adjunction between $E_*^R(E)$--comodules and $E_*$--modules
\[
\xymatrix{\Mod_{E_*}\ \ar@<0.5ex>[rr]^-{E_*^R(E)\otimes -} && \
\Comod_{E_*^R(E)} \ar@<0.5ex>[ll]^-{U},}
\]
where $U$ is the forgetful functor. Thus we obtain
\[
\Cohom^*_{E_*^R(E)}(E_*^R(Y), E_*^R(Z))\cong\Hom_{E_*}^*(E_*^R(Y),
Z_*),
\]
and the second statement follows from the first.
\end{proof}

We are now in a position to construct the $I$--adic tower. To
identify its homotopy and $F$--homology groups, however, we need a
few facts concerning universal coefficient and K\"unneth spectral
sequences. We recall these first.

Let $M$ and $N$ be $R$--modules. Then we have \cite[IV.4]{ekmm} a
universal coefficient spectral sequence of the form
\begin{equation}\label{ucss}
E_2^{*,*} = \Ext^{*,*}_{R^*}(M^*, N^*) \quad \Longrightarrow \quad
N^*_R(M).
\end{equation}
It is a conditionally convergent half-plane spectral sequence with
entering differentials, in the terminology of Boardman
\cite{boss}. If $M=N$, it is multiplicative with respect to
composition. It is shown in \cite[\S 7.3]{chr} that this spectral
sequences is in fact an Adams spectral sequence, with respect to
the ghost projective class in $\dr$.

Let us recall what we can say about the detection of a given
non-trivial map $f\: M\to N$. If the induced map $f_*\: M_*\to
N_*$ on homotopy groups is non-trivial, $f$ is represented in
\eqref{ucss} by $f_*$. This follows from the construction of the
spectral sequence. Assume now that $f_*$ is zero. Denoting the
homotopy cofibre of $f$ by $C$, we then have a short exact
sequence of $R_*$--modules
\begin{equation}\label{sest}
0 \lra N_* \lra C_* \lra M_{*-1} \lra 0.
\end{equation}
Let $t$ be the element of $\Ext^{1,-1}_{R_*}(M_*, N_*)$
representing it. We leave the verification of the following fact
to the reader.

\begin{prop}\label{detection}
If the element $t$ is non-trivial, it represents the map $f$ in
the $E_2$--term of \eqref{ucss}. Otherwise, $f$ has filtration
degree strictly greater than one.
\end{prop}

Let $P$ be a further $R$--module. There is \cite[IV.4]{ekmm} a
K\"unneth spectral sequence of the form
\begin{equation}\label{kss}
E^2_{*,*} = \Tor_{*,*}^{R_*}(P_*, M_*) \quad \Longrightarrow \quad
P_*^R(M).
\end{equation}
It is a strongly convergent half-plane spectral sequence with
exiting differentials.

We will need a statement about naturality. Namely, let $f\: M\to
N$ be a non-trivial map inducing zero on homotopy groups, as
above. Consider the K\"unneth spectral sequences \eqref{kss} and
\begin{equation}\label{kss2}
\bar E^2_{*,*} = \Tor_{*,*}^{R_*}(P_*, N_*) \quad \Longrightarrow
\quad P_*^R(N).
\end{equation}
As $f_*\: M_* \to N_*$ is zero, the map $f_*^r\: E^r_{*,*} \to
\bar E^r_{*,*}$ induced by $f$ on the spectral sequences is zero.
Because the spectral sequences converge strongly, this implies
that $P^R_*(f)\: P^R_*(M) \to P_*^R(N)$ induces zero on the
associated graded. Therefore, $P^R_*(f)$ induces a map of degree
$1$
\begin{equation}\label{gr+1}
\gr^*(P^R_*(M)) \lra \gr^{*+1}(P_*^R(N)).
\end{equation}
Similarly, the morphism $P^*_R(f)\: P^*_R(N) \to P^*_R(M)$ induces
a map
\begin{equation}\label{gr2+1}
\gr^*(P^*_R(N)) \lra \gr^{*+1}(P^*_R(M)).
\end{equation}
Indeed, to deduce that $P^*_R(f)$ induces zero on the associated
graded, it suffices to know that the respective spectral sequences
are conditionally convergent, see \cite[Lemma 5.6]{boss}. The
following statement is proved in a similar way as the Geometric
Boundary Theorem \cite[Thm.\@ 2.3.4]{ravenelgreen}.

\begin{prop}\label{gbt}
The homomorphism \eqref{gr+1} is represented on the K\"unneth
spectral sequences \eqref{kss} and \eqref{kss2} by a morphism
$E^r_{*,*} \to \bar E^r_{*-1,*}$ of degree $-1$. On $E^2$--terms,
it is given by the connecting homomorphism
\[
\Tor_{*,*}^{R_*}(P_*, M_*) \xrightarrow{\ \partial\ }
\Tor_{*-1,*}^{R_*}(P_*, N_*)
\]
associated to the short exact sequence \eqref{sest}. A similar
statement holds for the ho\-mo\-morphism \eqref{gr2+1} and the
universal coefficient spectral sequences converging to $P^*_R(N)$
and $P^*_R(M)$.
\end{prop}

Now consider the universal coefficient spectral sequence
\begin{equation}\label{ucssft}
E_2^{*,*} = \Ext_{R^*}^{*,*}(L^*, L^*) \ \Lra \ L^*_R(L).
\end{equation}
By Proposition \ref{exterior}, we have an isomorphism of
$L_*$--algebras
\[
\Ext^{*,*}_{R^*}(L^*, L^*) \cong \widehat\Lambda_{L^*}(f_0, f_1,
\ldots);
\]
the $f_j$ were explicitly defined there. We already know from
Proposition \ref{essential} that there is an isomorphism of
$L_*$--algebras
\[
L^*_R(L) \cong \widehat\Lambda_{L_*}(Q_0, Q_1, \ldots).
\]

\begin{prop}\label{repbocksteins}
The Bockstein operation $Q_j\: L \to \Sigma^{|x_j|+1} L$ is
represented in the spectral sequence \eqref{ucssft} by the element
$f_j$. Therefore, the spectral sequence collapses.
\end{prop}

\begin{proof}
Consider the cofibre sequence in $\dr$ associated to $Q_j$
\begin{equation}\label{cofbockstein}
L\xra{\ Q_j\ } \Sigma^{|x_j|+1} L \lra C_j \lra \Sigma L.
\end{equation}
By Proposition \ref{detection}, it suffices to verify two things.
Firstly, that $Q_j$ is trivial on homotopy groups. Secondly, that
the thereby obtained short exact sequence
\begin{equation}\label{sesqj}
0 \lra L_{*-|x_j|-1} \lra (C_j)_* \lra L_{*-1} \lra 0
\end{equation}
corresponds to $f_j$ in $\Ext_{R_*}^{1,*}(L_*, L_*)$. For the
first statement, note that the morphism induced by $Q_j$ on
homotopy groups is trivial for degree reasons.

For the second statement, recall that $Q_j$  is obtained by
smashing the composition
\begin{equation*}\label{decomp}
\bar\beta_j\: R/x_j \xrightarrow{\beta_j} \Sigma^{|x_j|+1} R
\xrightarrow{\rho_j} \Sigma^{|x_j|+1} R/x_j
\end{equation*}
with $\bigwedge_{i\neq j} R/x_i$. Applying the octahedral axiom
(see \cite[Appendix A]{hps} for instance) to the induced
decomposition of $Q_j$, we find that \eqref{sesqj} is isomorphic
to:
\begin{equation}\label{sesj}
0 \lra \Sigma^{|x_j|+1} L_* \xrightarrow{x_j} X_* = R_*/(x_j^2,
x_i;\ i\neq j) \lra L_* \lra 0
\end{equation}
Now the statement follows from the fact that the element
$f_j\in\Ext^{1,*}_{R_*}(L_*, L_*)$ corresponds to  \eqref{sesj}.
This verification is left to the reader.
\end{proof}

We are now ready to prove our main statement. Note the special
case $T=R$.

\begin{thm}\label{iadic}
Let $(R, T, F)$ be a regular triple with $F=T\wedge L$.

{\em (i)}\qua Applying $T\wedge -$ to the sequence
\eqref{resolutionsequence} yields a relative injective resolution
\begin{equation}\label{resolutionsequencet}
\ast \lra T \lra F \lra \Sigma I/I^2 \lra \Sigma^2 I^2/I^3 \lra
\cdots
\end{equation}
of $T$ with respect to the Adams injective class associated to $L$
in $\dr$ and hence with respect to the one associated to
$F=T\wedge L$.

{\em (ii)}\qua The Adams resolution associated to
\eqref{resolutionsequencet}, a diagram of $T$--module spectra, is
of the form:
\begin{equation}\label{iadictower}\begin{array}{c}
\xymatrix@!C=.8cm{ \cdots\ \ar[r]^{\iota} & I^3
\ar[r]^{\iota}\ar[d]_{\pi} & I^2 \ar[r]^{\iota}\ar[d]_{\pi} & I
\ar[r]^{\iota}\ar[d]_{\pi} & T\ar[d]_{\pi}
\\
& I^3/I^4\ar[lu]^{\epsilon^3}|-\circ &
I^2/I^3\ar[lu]^{\epsilon^2}|-\circ &
I/I^2\ar[lu]^{\epsilon^1}|-\circ & F \ar[lu]^{\epsilon^0}|-\circ}
\end{array}
\end{equation}
The induced diagram of homotopy groups is isomorphic to
\begin{equation}\label{htyiadictower}\begin{array}{c}
\xymatrix@!C=.8cm{ \cdots\ \ar[r]^{i} & I^3 \ar[r]^{i}\ar[d]_{p} &
I^2 \ar[r]^{i}\ar[d]_{p} & I \ar[r]^{i}\ar[d]_{p} & T \ar[d]_{p}
\\
& I^3/I^4 & I^2/I^3 & I/I^2\ & F }
\end{array}
\end{equation}
composed of short exact sequences. The induced diagram of
$F$--homology groups is isomorphic to the image of the diagram
\eqref{khomtower} under the functor $F_*^R(T)\otimes-$, equipped
with the total gradation.
\end{thm}

\begin{defn}\label{defiadic}
We call \eqref{iadictower} the {\em topological $I$--adic tower
over $T$}.
\end{defn}

\begin{proof}
(i)\qua First, we verify that \eqref{resolutionsequence} is a relative
injective resolution of $R$ with respect to $L$. Clearly, all
objects in the sequence \eqref{resolutionsequence}, except for
$R$, are $L$--injective. Let $N$ be an arbitrary $R$--module. By
definition of the Adams injective class associated to $F$ in
$\dr$, we have to check that the sequence
\begin{equation}\label{sequence}
\cdots \lra [L\wedge J^{s+1}/J^{s+2}, N]_R^{*-1}
\xrightarrow{(\wt\epsilon^s)^*} [L\wedge J^s/J^{s+1}, N]_R^* \lra
\cdots
\end{equation}
obtained by applying the functor $[L\wedge -, N]_R^*$ to
\eqref{resolutionsequence} is exact. The adjunction \eqref{adjrtf}
(for $F=L$) and Lemma \ref{kmodulemaps} imply
\[
[L\wedge J^s/J^{s+1}, N]_R^* \cong [J^s/J^{s+1}, L_R(L,N)]_R^*
\cong (L_*^R(J^s/J^{s+1}), N^{-*}_R(L))^*_{L_*}
\]
where $(-,-)^*_{L_*}=\Hom^*_{L_*}(-,-)$. But by Lemma
\ref{relinjstatement}, the $L$--homology of the sequence
\eqref{resolutionsequence} splits into short exact sequences of
$L_*$--modules. It follows that \eqref{sequence} is exact. Now
clearly \eqref{resolutionsequencet} is an $L$--resolution of $T$.
As $F$ is of the form $T\wedge L$, it is also an $F$--resolution.

(ii)\qua Consider the Adams resolution associated to the sequence
\eqref{resolutionsequence} first. It is of the form:
\begin{equation}\label{jadictower}\begin{array}{c}
\xymatrix@!C=.8cm{ \cdots\ \ar[r]^{\iota} & J^3
\ar[r]^{\iota}\ar[d]_{\pi} & J^2 \ar[r]^{\iota}\ar[d]_{\pi} & J
\ar[r]^{\iota}\ar[d]_{\pi} & R\ar[d]_{\pi}
\\
& J^3/J^4\ar[lu]^{\epsilon^3}|-\circ &
J^2/J^3\ar[lu]^{\epsilon^2}|-\circ &
J/J^2\ar[lu]^{\epsilon^1}|-\circ & L \ar[lu]^{\epsilon^0}|-\circ}
\end{array}
\end{equation}
We claim that the diagram of homotopy groups of the tower is
isomorphic to the diagram of homology groups of
\eqref{algebraiciadicover}, i.e.\ to:
\begin{equation}
\label{homalgebraiciadicover} \begin{array}{c} \xymatrix@!C=.8cm{
\cdots\ \ar[r]^{i} & J^3 \ar[r]^{i}\ar[d]^{p} & J^2
\ar[r]^{i}\ar[d]^{p} & J \ar[r]^{i}\ar[d]^{p} & R_*\ar[d]^{p} \\
& J^3/J^4 & J^2/J^3 & J/J^2 & L_*}
\end{array}
\end{equation}
It is clear that the sequence of homotopy groups of the triangle
\[
J \xra{\ i\ } R \xra{\ \pi\ } L \xra{\ \epsilon^0\ } \Sigma J
\]
splits into the short exact sequence
\[
\cale^0\:\quad  0 \lra J \xra{\ i\ } R_* \xra{\ p\ } L_* \lra 0.
\]
Consider the commutative diagram:
\[
\xymatrix@C=0.3cm@R=0.5cm{\Sigma^{-1}L \ar[rr]^{\tilde\epsilon^0}
\ar[rd]_{\epsilon^0} && J/J^2\\
& J\ar[ru]_{\pi_1}}
\]
Here, $\pi_1$ is the unique lift of $\tilde\epsilon^0$ to $J$. It
is unique because $J$ is the first term of the uniquely determined
Adams resolution associated to \eqref{resolutionsequence}. We need
to show that it induces the canonical projection $p_1$ on
coefficients. For this, we consider the pair of universal
coefficient spectral sequences:
\begin{equation*}\label{diagramss}\begin{array}{c}
\xymatrix@C=0cm@R=0.2cm{E_2^{*,*} & = & \Ext^{*,*}_{R_*}(J, J/J^2)
& \Longrightarrow & [J, J/J^2]_R^*
\\
\bar E_2^{*+1,*} & = & \Ext^{*+1,*}_{R^*}(L_*, J/J^2)&
\Longrightarrow & [L, J/J^2]^{*+1}_R}
\end{array}
\end{equation*}
By Proposition \ref{gbt}, the morphism $[J, J/J^2]^*_R \xra{\
(\epsilon^0)^*\ } [L, J/J^2]^{*+1}_R$ induced by $\epsilon^0$ is
represented by a morphism of degree one on the spectral sequences.
On $E_2$--terms it is given by the connecting homomorphism
\[
\epsilon^0\: \Ext^{*,*}_{R^*}(J, J/J^2) \lra
\Ext^{*+1,*}_{R^*}(L_*, J/J^2)
\]
associated to $\cale^0$. Let $x_1\in E_2^{*,*}$ and $y_0\in \bar
E_2^{*,*} $ be representatives of $\pi_1$ and $\tilde\epsilon^0$
respectively. We need to check that $x_1=p_1$. As
$\tilde\epsilon^0=(\epsilon^0)^*(\pi_1)$, we have $y_0 =
\epsilon^0(x_1)$ in case that $\epsilon^0(x_1)$ is non-trivial.
However, it follows from Theorem \ref{homtower} and Remark
\ref{cohomologydual} that $\epsilon^0$ is monomorphic. So it
suffices to identify $y_0$ and to show that $\epsilon^0(p_1)=y_0$.
Proposition \ref{repbocksteins} implies that
\[
\wt\epsilon^0=\bigvee_j Q_j \: L\to \Sigma J/J^2
\]
is represented by
\[
y_0 = \sum f_j\otimes \{x_j\}\in \Ext^{1,*}_{R_*}(L_*,
L_*)\otimes_{L^*} J/J^2\cong \Ext^{1,*}_{R_*}(L_*, J/J^2).
\]
The $f_j$ were defined in such a way that $y_0$ corresponds to the
short exact sequence
\begin{equation*}
\calf^1\: \quad  0 \lra  J/J^2 \xra{\ j_1\ }  R_*/ J^{2} \xra{\
q_{2}\ }  L_*  \lra 0.
\end{equation*}
It remains to verify that $\epsilon^0(p_1)$ is represented by
$\calf^1$. By definition, the connecting homomorphism $\epsilon^0$
maps the homomorphism $p_1\: J\to J/J^2$ to the pushout of the
short exact sequence $\cale^0$ along $p_1$. This pushout is
$\calf^1$, as the diagram of short exact sequences
\begin{equation*}\xymatrix@R=0.5cm{
\cale^0 : &    0 \ar[r] & J \ar[r]^{i_{1}}\ar[d]_{p_1} & R_*\ar[d]
\ar[r]^-{p_{0}} &
L_*\ar@{=}[d] \ar[r]& 0\\
\calf^1 : &   0 \ar[r] & J/J^{2} \ar[r]^{j_{1}} & R_*/ J^{2}
\ar[r]^-{q_{2}} & L_*  \ar[r] & 0}
\end{equation*}
shows. Here, the middle vertical map is the projection. We
conclude that
\begin{equation}\label{e1top}
J^2 \xrightarrow{\ \iota\ } J\xrightarrow{\ \pi\ } J/J^2
\xrightarrow{\ \epsilon^1\ } \Sigma J^2
\end{equation}
realizes the extension $\cale^1$. We can carry on by induction.

Let us now compute the homotopy groups of the Adams resolution
\eqref{iadictower} associated to \eqref{resolutionsequencet}. If
we smash \eqref{jadictower} with $T$, we obtain another tower
which qualifies as an Adams resolution associated to
\eqref{resolutionsequencet}. By uniqueness of Adams resolutions,
the two towers must be isomorphic. The product on $T_*$ and the
natural transformation $T_*\otimes_{R_*} -\to (T\wedge -)_*$
induce a commutative diagram
\[
\xymatrix@C=0.4cm{\cdots \ar[r]^-{0} & I^{s+1} \ar[r] & I^s \ar[r]
& I^s/I^{s+1} \ar[r]^-{0} & I^{s+1} \ar[r] & \cdots
\\
\cdots \ar[r]^-0 & T_*\otimes J^{s+1} \ar[r] \ar[d] \ar[u] &
T_*\otimes J^s \ar[r]\ar[d]  \ar[u] & T_*\otimes J^s/J^{s+1}
\ar[r]^-{0}\ar[d] \ar[u] & (T_*\otimes J^{s+1})_{*-1} \ar[r]\ar[d]
\ar[u] & \cdots
\\
\cdots \ar[r] & (T\wedge J^{s+1})_* \ar[r] & (T\wedge J^s)_*
\ar[r] & (T\wedge J^s/J^{s+1})_* \ar[r] & (T\wedge J^{s+1})_{*-1}
\ar[r] & \cdots}
\]
where $\otimes=\otimes_{R_*}$. We show that all vertical maps are
isomorphisms. For the upper row it suffices to note that
\[
\Tor_k^{R_*}(T_*, J^s/J^{s+1}) = \Tor_k^{R_*}(T_*,
R_*/J)\otimes_{L_*} J^s/J^{s+1} = 0
\]
for $k\neq 0$, which is a consequence of the assumption that the
sequence generating the ideal $J$ is regular on $T_*$. By
construction of $J^s/J^{s+1}$, we have
\[
(T\wedge J^s/J^{s+1})_* \cong I^s/I^{s+1}.
\]
Using this fact and starting with $(T\wedge J^0)_* = (T\wedge R)_*
\cong T_*$, we show inductively that $(T\wedge J^s)_* \cong I^s$,
using the Five-Lemma. We conclude that the homotopy groups are as
asserted.

The identification of the $L$--homology is easier. We show first
that $L_*^R(-)$ of the $J$-adic tower \eqref{jadictower} is
isomorphic to \eqref{khomtower}, endowed with the total gradation.
In the proof of Lemma \ref{relinjstatement} we have used the fact
that $L_*^R(-)$ of the sequence \eqref{resolutionsequence} is
isomorphic to the long exact sequence obtained from Theorem
\ref{homtower}
\[
0 \to L_* \xrightarrow{ p_* } HL_{*,*}(L_*) \xrightarrow{ \delta^0
} HL_{*-1,*}(J/J^2) \xrightarrow{ \delta^1 } HL_{*-2,*}(J^2/J^3)
\to \cdots
\]
endowed with the total gradation. It follows that we have
isomorphisms of $L_*$--modules:
\[
L^R_{*-1}(J) \cong \coker(\eta_*\: L_* \xra{\pi} L_*^R(L)) \cong
\coker(p_*) \cong HL_{*-1,*}(J)
\]
Furthermore, the morphism $(\pi_1)_*\: L^R_*(J)\to L^R_*(J/J^2)$
is injective, as it corresponds to $(p_1)_*\: HL_{*,*}(J) \to
HL_{*,*}(J/J^2)$. The fact that \eqref{e1top} is a cofibration
implies that $L_*^R(\Sigma J^2) \cong \coker(\pi_1)_*$. By an
inductive argument, the statement follows. In particular, we
obtain that all $L$--homology groups are $L_*$--free and hence
that all $F$--homology groups are $F_*$--free. Therefore
Proposition \ref{kunnethiso} implies that we have K\"unneth
isomorphisms
\[
F_*^R(T\wedge X)\cong F_*^R(T)\otimes F_*^R(X)
\]
for all terms $X$ in the $J$-adic tower. Thus we find that the
$F$--homology of the $I$-adic tower \eqref{iadictower} is as
claimed.
\end{proof}

\begin{rem}\label{adamsalgebraicresolution}
Consider the sequence in $\mathscr D_{R_*}$
\begin{equation}
\begin{array}{c}
\xymatrix{ T_* \ar[r] & F_* \ar[r]|-\circ & I/I^2 \ar[r]|-\circ &
I^2/I^3 \ar[r]|-\circ  & \ \cdots}
\end{array}
\end{equation}
derived from the algebraic $I$--adic tower \eqref{iadicintro2}. It
follows from Theorem \ref{homtower} and arguments analogous to the
ones in the preceding proof that it provides a relative injective
resolution of $T_*$ with respect to the Adams injective class
associated to $F_*$ (see Remark \ref{adamsalgebraic}).
Furthermore, \eqref{iadicintro2} is characterized (up to
isomorphism) as the Adams resolution associated to it.
\end{rem}

\begin{rem}\label{iadictowerunderT}
Having constructed an $I$--adic tower over $T$, we can easily
construct one under $T$, of the form:
\begin{equation}\label{quotiadic}\begin{array}{c}
\xymatrix{ T \quad \ar[r] & \quad\ldots\quad \ar[r]^-\rho & T/I^3
\ar[r]^-\rho & T/I^2 \ar[dl]^-{\theta^2} \ar[r]^-\rho & F
\ar[dl]^-{\theta^1}
\\
&& I^2/I^3 \ar[u]^-\nu & I/I^2 \ar[u]^-\nu}
\end{array}
\end{equation}
Namely, define spectra $T/I^s$ in $\dr$ as cofibres of the
canonical maps $I^s \to T$ and construct maps $T/I^{s+1}
\xrightarrow{\rho} T/I^s$ by making use of the octahedral axiom.
It is possible to choose these and compatible maps
\begin{equation} \label{holims}
\holim_s I^s \xrightarrow{\ \iota\ } T \xrightarrow{\ \rho\ }
\holim_s T/I^s \lra \holim_s I^s
\end{equation}
between the homotopy limits so that \eqref{holims} is a cofibre
sequence, see \cite[Rem.\@ after Prop.\@ 2.2.12]{hps}. The diagram
of homotopy groups of \eqref{quotiadic} is isomorphic to:
\begin{equation}\label{htyquotiadic}\begin{array}{c}
\xymatrix{ T_* \quad \ar[r] & \quad\ldots\quad \ar[r]^-\rho &
T_*/I^3 \ar[r]^-\rho & T_*/I^2  \ar[r]^-\rho & F_*
\\
&& I^2/I^3 \ar[u]^-\nu & I/I^2 \ar[u]^-\nu}
\end{array}
\end{equation}
Also, the diagram of $F$--homology groups of \eqref{quotiadic} is
isomorphic to $F_*^R(T)\otimes-$ of the one in Corollary
\ref{homquottower}, with the total gradation.
\end{rem}

\begin{notation}\label{that}
We write $I^\infty=\holim_s I^s$ and $\wh T=\holim_s T/I^s$ for
the homotopy limits of the $I$--adic towers. From Theorem
\ref{iadic}, we obtain $\wh{T}_*\cong \lim_s T_*/I^s =
(T_*)^\wedge_I$.
\end{notation}

We now aim to identify $\wh T$ under certain conditions as a
Bousfield localization with respect to $F$ in $\dr$ and in $\ds$.
First, we need to recall some terminology and basic facts about
Bousfield localizations.

Let $F$ be a module over a commutative $\ss$--algebra $R$. An
$R$--module $Z$ is called $F$--acyclic if $F\wedge Z\cong 0$. An
$R$--module $X$ is called $F$--local if $[Z, X]^*_R = 0$ for all
$F$--acyclic spectra $Z$. A map $f\: X \to Y$ of $R$--modules is
called an $F$--equivalence if $F\wedge f$ is an isomorphism. By
\cite[Thm.\@ 3.2.2]{hps}, we can associate to $F$ the Bousfield
localization and colocalization functors $L_F^R \: \dr \to \dr$
and $C_F^R\: \dr \to \dr$ respectively. We write $L_F$ and $C_F$
if $R=S$. The functors come with natural transformations $i^R\: 1
\to L_F^R$ and $q^R\: C_F^R \to 1$, where $1$ is the identity
functor. The pair $(L_F^R, i^R)$ is characterized by the fact that
$L_F^R X$ is $F$--local and that $i^R_X\: X \to L_F^R X$ is
initial among $F$--local objects under $X$. Also, $i^R_X$ is an
$F$--equivalence for each $X$. Similarly, $(C_F^R, q^R)$ is
characterized by $C_F^R X$ being $F$--acyclic and $q^R_X\: C_F^R X
\to X $ being terminal among $F$--acyclic objects over $X$. For
each $X$, there are natural cofibrations
\[
C_F^R X \xrightarrow{q^R_X} X \xrightarrow{i^R_X} L_F^R X \to
\Sigma C_F^R X.
\]
The pair $(L_F^R, i^R)$ determines $(C_F^R, q^R)$ by means of
these, see \cite[Lemma 3.1.6]{hps}.

The following fact is well-known. We include the proof, as it is
quite short.

\begin{lem}\label{localcor}
Let $X$ be an $F$--local in $\dr$. Then $X$ is $F$--local in
$\ds$.
\end{lem}

\begin{proof}
Let $Z$ be an $F$--acyclic spectrum. Then $F\wedge_R (R\wedge_\ss
Z) \cong F\wedge_\ss Z \cong 0$, so $R\wedge_\ss Z$ is
$F$--acyclic in $\dr$. Therefore, we have
\[
[Z, X]^* \cong [R\wedge_\ss Z, X]^*_R = 0
\]
which proves that $X$ is $F$--local in $\ds$.
\end{proof}

Recall that the ideal $I\lhd T_*$ is called invariant in
$T_*^R(T)$ iff $I\cdot T_*^R(T) = T_*^R(T)\cdot I$. Here we use
the left and the right actions of $T_*$ on $T_*^R(T)$ induced by
including $T$ via the unit $\eta_T\: R\to T$ as the right or the
left factor of $T\wedge_R T$ respectively. Similarly, $I$ is
defined to be invariant in $T_*(T)$ iff $I\cdot T_*(T)=T_*(T)\cdot
I$.

\begin{prop}\label{localization}
Assume that the sequence $S$ generating the ideal $I$ is finite.
Then the cofibre sequence \eqref{holims}
\[
I^\infty \xrightarrow{\ \iota\ } T \xrightarrow{\ p\ } \wh T \lra
\Sigma I^\infty
\]
can be identified with the triangle
\[
C_F^R T \xrightarrow{\ q_F^R\ } T \xrightarrow{\ i_F^R\ } L_F^R T
\lra \Sigma C_F^R T
\]
arising from localizing $T$ with respect to $F$ in $\dr$. If
moreover $I$ is invariant in $T_*(T)$, we can view \eqref{holims}
in $\ds$ as
\[
C_F T \xrightarrow{\ q_F\ } T \xrightarrow{\ i_F\ } L_F T \lra
\Sigma C_F T.
\]
\end{prop}

\begin{proof}
By what was said above, it suffices to exhibit $p\: T\to \wh T$ as
the Bousfield localization of $T$ in $\dr$ with respect to $F$ to
prove the first statement. By construction, all the objects
$T/I^s$ lie in the thick subcategory of $\dr$ generated by $F$ and
are therefore $F$--local. The class of $F$--local objects is
closed under homotopy limits, hence $\wh T$ is $F$--local. To
prove that $p$ is an $F$--equivalence, we first show that the
composition
\begin{equation}\label{firststep}
F\wedge T \lra \holim_s(F\wedge T/I^s)
\end{equation}
of $F\wedge p\: F\wedge T \to F\wedge \wh T$ with the natural map
\begin{equation}\label{can}
F\wedge\wh T = F\wedge\holim_s(T/I^s) \lra \holim_s(F\wedge T/I^s)
\end{equation}
is an equivalence. Consider the Milnor type short exact sequence
\[
0 \lra \lim\nolimits^1_s F^R_{*+1}(T/I^s) \lra (\holim_s F\wedge
T/I^s)_* \lra \lim\nolimits_s F^R_*(T/I^s) \lra 0
\]
of $F_*$--modules, see \cite[Prop.\@ 2.2.11]{hps}. The
identification of the $F$--homology of the tower under $T$ in
Remark \ref{iadictowerunderT} together with Corollary
\ref{homquottower} imply that
\[
\lim\nolimits_s F^R_*(T/I^s) \cong F^R_*(T)
\]
and $\lim^1_s F^R_*(T/I^s)=0$. Therefore \eqref{firststep} is an
equivalence. It remains to prove that \eqref{can} is one as well.
Choose a cell approximation $\Gamma F$ of $F$ and a filtration
$(\Gamma F)_k$ by finite subcomplexes \cite[Lemma III.2.3]{ekmm}.
Fix a $K$ such that
\[
x_i\: ((\Gamma F)_k)_* \lra ((\Gamma F)_k)_*
\]
is trivial for $0\leq i\leq n$ for all $k>K$. Let $k>K$. It
follows by induction over the elements $x_i$ of the sequence
$S=(x_0, \ldots, x_n)$ that
\begin{equation}\label{ff}
((\Gamma F)_k)_*^R(F)_* \cong ((\Gamma F)_k)_*^R(T) \otimes_{L_*}
\Lambda_{L_*}(a_0', \ldots, a_n')
\end{equation}
with $|a'_i|=|x_i|+1$. Consider the map
\[
((\Gamma F)_k)^R_*(T\wedge Q_j)\: ((\Gamma F)_k)^R_*(F) \lra
((\Gamma F)_k)^R_{*-|x_j|-1}(F)
\]
induced by the Bockstein $Q_j\: L \to \Sigma^{|x_j|+1} L$. Under
\eqref{ff}, it corresponds to
\[
M_{k*}\otimes\frac{\partial}{\partial a'_j}\: M_{k*} \otimes_{L_*}
\Lambda_{L_*}(a_0', \ldots, a_n') \lra M_{k*} \otimes_{L_*}
\Lambda_{L_*}(a_0', \ldots, a_n')
\]
where $M_{k*}=((\Gamma F)_k)_*^R(T)$. It follows that we can
identify the images of the sequence \eqref{resolutionsequence}
under the two functors $((\Gamma F)_k)_*(T\wedge-)$ and $M_{k*}
\otimes_{L_*} L_*^R(-)$. This implies that $((\Gamma F)_k)_*(-)$
of the $I$--adic tower over $T$ is isomorphic to $M_{k*}
\otimes_{L_*} L_*^R(-)$ of the Adams resolution associated to
\eqref{resolutionsequence}. In particular, we find that
\[
\lim\nolimits_s ((\Gamma F)_k)^R_*(T/I^s) \cong M_{k*}, \
\lim\nolimits^1_s ((\Gamma F)_k)^R_*(T/I^s)=0.
\]
Thus, a Milnor type short exact sequence implies that
\[
(\holim_s (\Gamma F)_k \wedge T/I^s)_* \cong M_{k*}.
\]
As $(\Gamma F)_k$ is finite, it follows that $((\Gamma F)_k\wedge
\wh T)_* \cong M_{k*}$. Therefore, we have
\[
F^R_*(\wh T) \cong \colim_k((\Gamma F)_k\wedge \wh T)_* \cong
\colim_k M_{k*} \cong F^R_*(T).
\]
Now assume that $I$ is invariant in $T_*(T)$. By Lemma
\ref{localcor}, $\wh T$ is $F$--local in $\ds$. We need to show
that $F_*(T) \to F_*(\wh{T})$ is an isomorphism. Similar to
before, we first prove that
\begin{equation}\label{anotherstep}
(\holim_s F\wedge T/I^s)_* \cong F_*(T)
\end{equation}
as follows. As $I$ is invariant, the right $T_*$--action on
$F_*(T)$ factors through $F_*$. By induction over the elements
$x_i$ we obtain
\begin{equation}\label{ff'}
F_*(F)\cong F_*(T)\otimes\Lambda_{F_*}(a''_0,\ldots, a''_n)
\end{equation}
with $|a''_i|=|x_i|+1$. Analogous to the above, we find that
\[
F_*(T\wedge Q_j)\: F_*(F) \lra F_{*-|x_j|-1}(F)
\]
corresponds under \eqref{ff'} to:
\[
F_*(T)\otimes \frac{\partial}{\partial a''_j}\:
F_*(T)\otimes\Lambda_{F_*}(a''_0,\ldots, a''_n) \lra
F_*(T)\otimes\Lambda_{F_*}(a''_0,\ldots, a''_n)
\]
Arguing along the same lines as before, this implies
\eqref{anotherstep}. To prove that
\[
F_*(\wh T) = F_*(\holim_s T/I^s) \lra (\holim_s F\wedge_\ss
T/I^s)_*
\]
is an isomorphism, we take the sequence $(\Gamma F)_k$ from above,
pick a $K'$ such that
\[
((\Gamma F)_k \wedge x_i)\: ((\Gamma F)_k \wedge_\ss T)_* \lra
((\Gamma F)_k \wedge_\ss T)_*
\]
is trivial for $0\leq i\leq n$ for all $k>K'$ and argue similarly
as above.
\end{proof}

The ideal $I_n$ is invariant in $BP_*(BP)$ \cite{ravenelgreen}. By
Landweber exactness of $E(n)_*$, we have
\[
E(n)_*(E(n))\cong E(n)_*\otimes_{BP_*}BP_*(BP)\otimes_{BP_*}
E(n)_*
\]
and hence $I_n$ is also invariant in $E(n)_*(E(n))$. Thus we
obtain:

\begin{cor}
The natural maps $E(n)\to\wh{E(n)}$ and $BP\to \wh{BP}$, obtained
by considering the regular triples $(MU, E(n), K(n))$ and $(MU,
BP, P(n))$, are Bousfield localizations with respect to $K(n)$ and
$P(n)$ respectively.
\end{cor}

\begin{rem}
As $P(n)$ is connective, the second localization can also be
interpreted as the $p$-completion of $BP$ \cite[Thm.\@
3.1]{bousfield}.
\end{rem}

We record \eqref{ff} from the proof of Proposition
\ref{localization} and a corresponding statement for cohomology:

\begin{lem}\label{ffhc}
If $I$ is invariant in $T_*(T)$ and $S=(x_0, \ldots, x_n)$ is
finite, we have an isomorphism
\[
F_*(F) \cong F_*(T)\otimes \Lambda_{F_*}(a_0', \ldots, a_n')
\]
of $F_*$--modules, with $|a_i'|=|x_i|+1$. If furthermore $F_*(T)$
is $F_*$--free, there is an isomorphism of $F^*$--modules
\[
F^*(F) \cong F^*(T)\otimes \Lambda_{F^*}(Q_0, \ldots, Q_n).
\]
In this case, the exterior algebra $F^*\otimes L^*_R(L)$ maps
isomorphically onto the subalgebra $\Lambda_{F^*}(Q_0, \ldots,
Q_n)$ under the natural map of $F^*$--algebras
\[
F^*\otimes L^*_R(L) \to F^*_R(F) \to F^*(F).
\]
\end{lem}

\begin{rem}\label{iadicallycomplete}
Because of $\wh T_* \cong (T_*)^\wedge_I$, the map $T\to \wh T$ is
an equivalence in case that $T_*$ is $I$--adically complete. We
have seen in the proof of Proposition \ref{localization} that $\wh
T$ is $F$--local in $\dr$ and hence in $\ds$ by Lemma
\ref{localcor}. So we obtain:
\end{rem}

\begin{prop}
Assume that $(R, T, F)$ is a regular triple, with $F_*\cong
T_*/I$. If $T_*$ is $I$--adically complete, it follows that $T$ is
local with respect to $F$, in $\dr$ as well as in $\ds$.
\end{prop}

In \cite{bwart}, the notation $\wh{E(n)}$ denotes the spectrum
representing completed John\-son--Wilson theory. Recall that the
latter is defined on finite spectra $X$ by:
\[
\wh{E(n)}^*(X) = \lim_s E(n)^*(X)/(I_n^s\. E(n)^*(X)) \cong
\lim_s(E(n)^*/I_n^s)\otimes_{E(n)^*} E(n)^*(X)
\]
This determines the theory uniquely, because its coefficients are
linearly compact with respect to the $I_n$--adic topology.

\begin{prop}\label{enhat}
The homotopy limit $\wh{E(n)}$ of the $I_n$--adic tower under
$E(n)$ is isomorphic in $\ds$ to the spectrum representing
completed Johnson--Wilson theory.
\end{prop}

\begin{proof}
We claim that it suffices to construct natural morphisms
\begin{equation}\label{nattranscompl}
E(n)^*/I_n^s\otimes_{E(n)^*} E(n)^*(X) \lra (E(n)/I_n^s)^*(X)
\end{equation}
for finite complexes $X$ and verify that they are isomorphisms on
suspensions of spheres. Namely, we may then take inverse limits
over $s$ on both sides. On the left hand side, we obtain completed
Johnson--Wilson theory. For the right hand side, we have a Milnor
type short exact sequence \cite[Prop.\ 2.2.11]{hps}:
\[
0 \lra \lim\nolimits^1_s (E(n)/I_n^s)^{*-1}(X) \lra \wh{E(n)}^*(X)
\lra \lim\nolimits_s (E(n)/I_n^s)^*(X) \lra 0
\]
We claim that the $\lim^1$--term vanishes. By \cite[Thm.\@
7.1]{jensen}, it suffices to show that the $E(n)^*$--modules
\[
M_s^*=(E(n)/I_n^s)^*(X)
\]
are linearly compact with respect to the $I_n$--adic topology. By
Lemma \ref{powersofideals} below, the $E(n)^*$--action on $M_s^*$
factors through $E(n)^*/I_n^s$. So $M_s^*$ is discrete and
therefore Hausdorff, as $M_s^*$ is finitely generated. With
\cite[\S 7, Prop.\@ B]{jensen}, it follows that $M_s^*$ is
linearly compact, because $E(n)^*/I_n^s$ is so. Thus, completed
Johnson--Wilson theory and $\wh{E(n)}^*(-)$ coincide on finite
complexes. As the coefficients are linearly compact, this implies
the claim.

The existence of morphisms as in \eqref{nattranscompl} is
guaranteed by Lemma \ref{powersofideals} below.
\end{proof}

\begin{lem}\label{powersofideals}
For any $R$--module $M$, the ideal $I^s\lhd T_*$ lies in the
annihilator of $(T/I^s)_R^*(M)$.
\end{lem}

The proof will be given in the next section.

\section{Higher Bockstein spectral sequences}\label{sectionhbss}

\subsection{Setting up the spectral sequence; convergence}

Let $(R, T, F)$ be a regular triple, with $F=T\wedge L$. In the
last section, we have constructed the $I$--adic tower as an Adams
resolution of $T$ with respect to $F$ in $\dr$. Let $M$ be an
$R$--module. If we apply $-\wedge M$ to the tower, we obtain an
Adams resolution of $T\wedge M$. It is of the form:
\begin{equation}\label{adamsresolutiongeneral}
\begin{array}{c}
\xymatrix@!=.5cm{ T\wedge M\ar[rd] && I\wedge\ar[rd] M\ar[ll]&&
I^2\wedge M\ar[rd] \ar[ll]&&I^3\wedge M\ar[ll] &
\ar@{}[d]|{.\,.\,.}
\\
& F\wedge M\ar[ru]|-\circ && I/I^2\wedge M \ar[ru]|-\circ &&
I^2/I^3\wedge M \ar[ru]|-\circ& & }
\end{array}
\end{equation}
Applying $[Y,-]^*_R$, where $Y$ is an $R$--module, yields an exact
couple, hence a half plane spectral sequence with entering
differentials \cite[\S 7]{boss}.

Let us first recall what we can say about convergence. In general,
the homotopy limit $\holim_s I^s\wedge M$ need not be trivial, so
we do not have a good convergence behaviour with respect to the
target $[Y, T\wedge M]_R^*$. Boardman shows in \cite[\S 15]{boss}
how this can be remedied. Namely, define $R$--modules $Z^s(M)$ by
forming cofibre sequences
\begin{equation}\label{zsm}
\holim_t (I^t \wedge M) \lra I^s\wedge M \lra Z^s(M)
\end{equation}
where the first map is the canonical one. With \cite[Rem.\@ after
Prop.\@ 2.2.12]{hps}, we find that
\[
Z^0(M) \cong \holim_s (T/I^s\wedge M).
\]
Now construct compatible maps between the $Z^s(M)$ via the
octahedral axiom. The maps in \eqref{zsm} then induce a morphism
of \eqref{adamsresolutiongeneral} into the tower
\begin{equation}\label{adamsresolutionreplace}
\begin{array}{c}
\xymatrix@!=.5cm{ Z^0(M)  \ar[rd] &&  Z^1(M) \ar[rd] \ar[ll] &&
Z^2(M) \ar[rd] \ar[ll] && Z^3(M) \ar[ll] & \ar@{}[d]|{.\,.\,.}
\\
& T \wedge M\ar[ru]|-\circ && I/I^2\wedge M \ar[ru]|-\circ &&
I^2/I^3\wedge M \ar[ru]|-\circ& & }
\end{array}
\end{equation}
It induces an isomorphism of the respective spectral sequences
obtained after applying $[Y,-]_R^*$, because the $E_1$--terms
coincide. In the second tower, we have $\holim_s Z^s(M)=0$. Hence
the spectral sequence converges conditionally to
\[
[Y, Z^0(M)]^*_R \cong [Y, \holim_s (T/I^s\wedge M)]_R^*.
\]
The homotopy limit $\holim_s (T/I^s\wedge M)$ can be abstractly
characterized as the {\em $F$--nilpotent completion\/}
$F_R^\wedge(T\wedge M)$ of $T\wedge M$ in $\dr$. By this, we mean
the obvious generalization to $\dr$ of nilpotent completion with
respect to a ring spectrum in $\ds$, as defined in
\cite{bousfield}. In particular, $F_R^\wedge(T\wedge M)$ is
$F$--local for each $M$. This can be seen directly as follows. By
induction $T/I^s$ and hence $T/I^s\wedge M$ are $F$--nilpotent,
\ie{} lie in the ideal of $\dr$ generated by $F$, in the
terminology of \cite{hps}. In particular, they are $F$--local.
Therefore $F_R^\wedge(T\wedge M)$ is $F$--local, as a homotopy
limit of $F$--local $R$--modules. It follows that the canonical
map
\[
T\wedge M\to F_R^\wedge(T\wedge M)
\]
factors uniquely into the $F$--localization map $i^R_F\: T\wedge M
\to L_F^R(T\wedge M)$ composed with a map
\[
\varphi_{T\wedge M}\: L_F^R(T\wedge M) \lra F_R^\wedge(T\wedge M).
\]
For each finite $R$--module $M$, the natural map $\wh T\wedge M
\to F_R^\wedge(T\wedge M)$ is an equivalence, as an induction over
the cells of $M$ shows.

If the sequence $S$ generating $I$ is finite, the map $\varphi_R$
is an equivalence by Proposition \ref{localization}. If $T=R$, the
map
\[
\varphi_{R\wedge M}\:  L_F^R(T\wedge M) \cong L_F^R(M) \lra
F_R^\wedge(M)\cong F_R^\wedge(T\wedge M)
\]
is an equivalence for all $M$. Namely, we have K\"unneth
isomorphisms
\[
F_*^R(R/I^s\wedge M) \cong F_*^R(R/I^s)\otimes F_*^R(M)
\]
for all $s$, because $F_*^R(R/I^s)$ is $F_*$--free. Now we have
seen in the proof of Proposition \ref{localization} that $\lim_s
F_*^R(R/I^s)\cong F_*$ and $\lim_s^1 F_*^R(R/I^s)=0$. Therefore,
we obtain
\[
\lim\nolimits_s F_*^R(R/I^s\wedge M) \cong F_*^R(M), \
\lim\nolimits^1_s F_*^R(R/I^s\wedge M)=0.
\]
Using a Milnor short exact sequence and the fact that $F$ is a
finite $R$--module, the claim follows.

\subsection{The $E_2$--term; bicomodules}

If $F_*^R(Y)$ is $F_*$--free, it follows from Proposition
\ref{kmodulemaps} that the $E_1$--term of the spectral sequence is
given by:
\begin{equation}\label{e1adams}
E_1^{s,*} = [Y, I^s/I^{s+1}\wedge M]^*_R \cong
\Hom^*_{F_*}(F_*^R(Y), (I^s/I^{s+1}\wedge M)_*)
\end{equation}
If $L$ (and hence $F$) is commutative as an $R$--ring spectrum and
$F_*^R(T)$ (and hence $F_*^R(F)\cong F_*^R(T)\otimes F_*^R(L)$) is
$F_*$--flat, the same proposition yields
\begin{align*}
E_1^{s,*} & \cong \Cohom^*_{F_*^R(F)}(F_*^R(Y),
F_*^R(I^s/I^{s+1}\wedge M))).
\end{align*}
To simplify this expression, we need to gain a better
understanding of the $F_*^R(F)$--comodule structure. We consider
the following general setup.

Let $A$ and $B$ be coalgebras over some commutative ring $R$ with
coactions $\Delta_A$ and $\Delta_B$ and counits $\epsilon_A$ and
$\epsilon_B$. We write $\otimes$ for $\otimes_R$ in the following.
Let $\tau\: A\otimes B \to B\otimes A$ be a given symmetry
isomorphism of $R$--modules with $\tau^2=1$. Then we define a left
$(A,B)$--bicomodule with respect to $\tau$ to be an $R$--module
$M$ which is a left comodule both over $A$ and $B$ in such a way
that the diagram
\[
\xymatrix@C=0.2cm{M \ar[rr]^-{\gamma_A}\ar[d]^-{\gamma_B} &&
A\otimes M\ar[d]_-{A \otimes \gamma_B}
\\
B\otimes M \ar[rd]^-{B\otimes \gamma_A} && A\otimes B\otimes
M\ar[ld]_-{\tau\otimes M}
\\
& B\otimes A\otimes M}
\]
commutes. Here $\gamma_A$ and $\gamma_B$ denote the given
coactions. If we endow $A\otimes B$ with the coproduct
\[
\Delta\: A\otimes B \xra{\ \Delta_A\otimes\Delta_B\ } A^{\otimes
2}\otimes B^{\otimes 2}\xra{\ A\otimes\tau \otimes B\ } (A\otimes
B)^{\otimes 2},
\]
we find that
\[
M\xra{\gamma_A} A\otimes M \xra{A\otimes \gamma_B} A\otimes
B\otimes M
\]
defines a coaction of $A\otimes B$ on $M$. To check the
coassociativity axiom, we need to verify that the two paths from
the top left to the bottom right corner along the outer edges of
the following diagram yield the same map:
\[
\def\objectstyle{\scriptstyle} \xymatrix@C=1.8cm{ M
\ar@{}[dr]|*+[o][F-]{1} \ar[r]^-{\gamma_A}\ar[d]_-{\gamma_A} &
A\otimes M \ar@{}[drr]|*+[o][F-]{2} \ar[rr]^-{A\otimes\gamma_B}
\ar[d]_-{A\otimes\gamma_A} && A\otimes B\otimes M
\ar[d]_-{A\otimes B\otimes\gamma_A}
\\
A\otimes M \ar@<-1ex>@{}[dr]|*+[o][F-]{3} \ar[r]^-{\Delta_A\otimes
M} \ar[d]_-{A\otimes \gamma_B} & A^{\otimes 2} \otimes M
\ar@<-1ex>@{}[dr]|*+[o][F-]{4} \ar[r]^-{A^{\otimes 2}
\otimes\gamma_B} \ar[d]_-{A^{\otimes 2} \otimes \gamma_B} &
A^{\otimes 2} \otimes B\otimes M \ar@<-1ex>@{}[dr]|*+[o][F-]{5}
\ar[r]^-{A\otimes \tau\otimes M} \ar[d]_-{A^{\otimes 2}\otimes
B\otimes \gamma_B} & A\otimes B\otimes A\otimes M
\ar[d]_-{A\otimes B\otimes A\otimes\gamma_B}
\\
A\otimes B\otimes M \ar[r]^-{\Delta_A\otimes B\otimes M} &
A^{\otimes 2} \otimes B\otimes M \ar[r]^-{A^{\otimes 2}
\otimes\Delta_B \otimes M} & A^{\otimes 2} \otimes B^{\otimes 2}
\otimes M \ar[r]^-{A \otimes \tau\otimes B\otimes M} & (A\otimes
B)^{\otimes 2} \otimes M}
\]
This is the case, as the diagram is built from commuting squares:
Squares \textcircled{1} and \textcircled{4} commute because
$\gamma_A$ and $\gamma_B$ are coactions; square \textcircled{2}
commutes because $M$ is an $(A,B)$-bicomodule with respect to
$\tau$; the remaining two squares commute trivially. The
verification of the counit axiom is easy.

Vice versa, an $A\otimes B$--comodule is an $(A, B)$--comodule
with respect to $\tau$. Namely, the map
\[
\gamma_A\: M \xra{\gamma} A\otimes B \otimes M
\xra{A\otimes\epsilon_B \otimes M} A\otimes M,
\]
where $\gamma$ denotes the given $A\otimes B$--coaction, defines a
coaction of $A$ on $M$. Similarly, we define a coaction $\gamma_B$
of $B$. It is easy to check that the coassociativity and the
counit axioms are satified. To show that $M$ is an $(A,
B)$--bicomodule with respect to $\tau$, we check that both
\begin{equation}\label{comp1}
M\xra{\gamma_A} A\otimes M \xra{A\otimes\gamma_B} A\otimes B
\otimes M
\end{equation}
and
\begin{equation}\label{comp2}
M \xra{\gamma_B} B\otimes M \xra{B\otimes\gamma_A} B\otimes
A\otimes M \xra{\tau\otimes M} A\otimes B\otimes M
\end{equation}
coincide with the coaction $\gamma$. For \eqref{comp1}, this can
be seen by considering the commutative diagram:
\[
\xymatrix{M \ar[r]^-{\gamma}\ar[d]^-\gamma & A\otimes B\otimes M
\ar[r]^-{\epsilon_A \otimes B\otimes M} \ar[d]^-{A\otimes B\otimes
\gamma} & A\otimes M \ar[d]^-{A\otimes\gamma}
\\
A\otimes B\otimes M \ar[r]^-{\Delta\otimes M} \ar[rd]^-{=} &
(A\otimes B)^{\otimes 2}\otimes
M\ar[d]^-{A\otimes\epsilon_B\otimes \epsilon_A \otimes B\otimes M}
& A\otimes A\otimes B\otimes M \ar[d]^-{A\otimes\epsilon_A\otimes
B\otimes M}
\\
& A\otimes B\otimes M \ar[r]^-{=} & A\otimes B\otimes M }
\]
A similar diagram shows that \eqref{comp2} is $\gamma$. From now
on, we assume that $\tau$ is fixed and don't distinguish between
$A\otimes B$--comodules and $(A,B)$--bicomodules with respect to
$\tau$.

If $N$ is an $A$--comodule with coaction $\gamma_A$, then the
composition
\[
B\otimes N \xra{B\otimes\gamma_A} B\otimes A\otimes N
\xra{\tau\otimes N} A\otimes B\otimes N
\]
defines an $A$--coaction on $B\otimes N$. Together with the
coaction of $B$ given by
\[
B\otimes N \xra{\Delta_B\otimes N} B\otimes B\otimes N,
\]
$B\otimes N$ can be checked to be an $A\otimes B$--comodule. We
say that $B\otimes N$ is obtained by extending the $A$--coaction
of $N$ to $A\otimes B$. In a similar way, we can define an
extension $A\otimes N'$ of a given coaction of $B$ on $N'$ to
$A\otimes B$. It is then easy to show that we have the following
adjunctions for an $A\otimes B$--comodule $M$:
\begin{align*}
\Cohom_{A}(M, N) & \cong \Cohom_{A\otimes B}(M, B\otimes N)
\\
\Cohom_{B}(M, N')& \cong \Cohom_{A\otimes B}(M, A\otimes N')
\end{align*}

\begin{prop}
Let $(R, T, F)$ be a regular triple with $F=T\wedge L$. Assume
that $L$ (and hence $F$) is commutative and that $F_*^R(T)$ (and
hence $F_*^R(F))$ is $F_*$--flat. Then there are natural morphisms:
\begin{align}
\label{ft} F_*^R(M) & \lra F_*^R(T) \otimes F_*^R(M)\\
\label{fl} F_*^R(M) & \lra F_*^R(L) \otimes F_*^R(M)
\end{align}
For $M=T$ and $M=L$, they induce coalgebra structures on
$F_*^R(T)$ and $F_*^R(L)$ respectively. With respect to these,
\eqref{ft} and \eqref{fl} define natural coactions on $F_*^R(M)$.
Furthermore, there is a symmetry isomorphism
\[
\tau\: F_*^R(T)\otimes F_*^R(L) \cong F_*^R(L)\otimes F_*^R(T)
\]
with $\tau^2=1$ and with respect to which $F_*^R(M)$ is an
$(F_*^R(T), F_*^R(L))$--bicomo\-dule. The corresponding
$F_*^R(T)\otimes F_*^R(L)$--coaction coincides under the
K\"un\-neth isomorphism $F_*^R(F)\cong F_*^R(T)\otimes F_*^R(L)$
with the natural $F_*^R(F)$--coaction.
\end{prop}

We content ourselves with defining the natural maps \eqref{ft},
\eqref{fl} and the symmetry $\tau$. The verification of the
statements in the proposition is a laborious but straightforward
task.

The map \eqref{ft} is defined as the composition
\[
F_*^R(M) \xra{F_*^R(\eta_T\wedge M)} F_*^R(T\wedge M) \cong
F_*^R(T)\otimes F_*^R(M)
\]
where $\eta_T\: R\to T$ is the unit of $T$ and the isomorphism is
the inverse of:
\begin{equation}\label{twistedkunnethft}
F_*^R(T)\otimes F_*^R(M) \to (T_0\wedge L_0\wedge T_1 \wedge
T_2\wedge L_1\wedge M)_* \to (T_0 \wedge L_{01} \wedge T_{12}
\wedge M)_*
\end{equation}
Here the first map is the canonical one. The purpose of the
indices of the various copies of $T$ and $L$ is to indicate in
what way these are multiplied under the second map. The right
$F_*$--action on $F_*^R(T)$ used to form the tensor product is
defined as follows. An element $\gamma_1\wedge\gamma_2\in (T\wedge
L)_*=F_*$ acts on $x\in F_*^R(T)$ as
\begin{multline*}
R\xra{x} F\wedge T  \cong (T\wedge L\wedge R)\wedge T\wedge R
\xra{T\wedge L\wedge \gamma_2\wedge T \wedge \gamma_1}\\
(T\wedge L\wedge L)\wedge T\wedge T \xra{T\wedge\mu_L\wedge \mu_T}
(T\wedge L)\wedge T = F\wedge T
\end{multline*}
where $\mu_L$ and $\mu_T$ are the products on $L$ and $T$
respectively. Similarly, we define \eqref{fl} as
\[
F_*^R(M)\xra{F_*^R(\eta_L \wedge M)} F_*^R(L)\otimes F_*^R(M)
\cong F_*^R(L\wedge M),
\]
where $\eta_L\: R\to L$ is the unit of $L$ and where the
isomorphism is the inverse of:
\begin{equation}\label{twistedkunnethfl}
F_*^R(L)\otimes F_*^R(M) \to (T_0\wedge L_0\wedge L_1 \wedge
T_1\wedge L_2\wedge M)_* \to (T_{01} \wedge L_0 \wedge L_{12}
\wedge M)_*
\end{equation}
The symmetry $\tau$ is a composition of isomorphisms:
\[
F_*^R(T)\otimes F_*^R(L) \cong F_*^R(T\wedge L) \cong
F_*^R(L\wedge T) \cong F_*^R(L)\otimes F_*^R(T)
\]
The first map is \eqref{twistedkunnethft} for $M=L$, the second is
induced by the switch $T\wedge L\cong L\wedge T$ and the third is
the inverse of \eqref{twistedkunnethfl} for $M=T$.

\begin{rem}\label{teinfty}
If $T$ is a {\em commutative} $\ss$--algebra, $F=T\wedge L$ is
automatically a $T$--ring spectrum. It is commutative as such if
$L$ is commutative as an $R$--ring spectrum. In this case,
$F_*^T(-)$ defines a multiplicative homology theory on $\dt$
taking values in $F_*^T(F)$--comodules. It can be checked that the
coaction of $F_*^T(F) \cong F_*^R(L)$ on $F_*^T(T\wedge M) \cong
F_*^R(M)$ defined in this way is the same as the one considered
above.
\end{rem}

We leave the verification of the following fact to the reader as
well:

\begin{prop}
For any $R$--module $N$, we have an isomorphism of
$F_*^R(F)$--comodules
\[
F_*^R(T\wedge N) \cong F_*^R(T) \otimes F_*^R(N),
\]
where the coaction on the right hand side is the one obtained by
extending the $F_*^R(L)$--coaction on $F_*^R(N)$ to $F_*^R(F)$.
\end{prop}

Let us come back to the identification of the $E_1$--term of our
spectral sequence. It follows from the definition of
$I^s/I^{s+1}$, the preceding discussion and the fact that
$F_*^R(J^s/J^{s+1})$ is $F_*$--free that we have isomorphisms:
\begin{align*}
E_1^{s,*} & \cong \Cohom^*_{F_*^R(F)}(F_*^R(Y),
F_*^R(I^s/I^{s+1}\wedge M)))\\
& = \Cohom^*_{F_*^R(F)}(F_*^R(Y),
F_*^R(T\wedge J^s/J^{s+1}\wedge M)))\\
& \cong \Cohom^*_{F_*^R(L)}(F_*^R(Y),
F_*^R(J^s/J^{s+1}\wedge M))\\
& \cong \Cohom^*_{F_*^R(L)}(F_*^R(Y), F_*^R(J^s/J^{s+1}) \otimes
F_*^R(M))
\end{align*}
Thus the $E_2$--term is the cohomology of the complex obtained by
applying the functor $\Cohom^*_{F_*^R(L)}(F_*^R(Y),-)$ to the
sequence:
\begin{equation}\label{flcomodules}
0 \to F_*^R(M) \to F_*^R(L)\otimes F_*^R(M) \to F_*^R(\Sigma
J/J^2)\otimes F_*^R(M) \to \cdots
\end{equation}
We claim that this implies that
\[
E_2^{*,*}\cong\Coext^{*,*}_{F_*^R(L)}(F_*^R(Y),F_*^R(M)).
\]
To see this, note that \eqref{flcomodules} is the sequence of
homotopy groups of
\begin{equation}\label{fm}
\ast \lra F\wedge M \lra F\wedge L\wedge M \lra F\wedge \Sigma
J/J^2 \wedge M \lra \cdots
\end{equation}
obtained from \eqref{resolutionsequence} by applying
$F\wedge-\wedge M$. We have noted in the proof of Theorem
\eqref{iadic} that \eqref{resolutionsequence} is a relative
injective resolution with respect to $L$. Therefore, the sequence
\eqref{fm} is split, \ie{} exact with respect to the split
injective class in $\dr$. It follows that the sequence
\eqref{flcomodules} is split over $F_*$. This implies that it is a
relative injective resolution of $F_*^R(M)$ over $F_*^R(L)$ and
thus, as $F_*^R(Y)$ is $F_*$--projective, the $E_2$--term is as
claimed.

\subsection{The Higher Bockstein spectral sequence}

Let now $X$ be a spectrum, \ie{} an object of $\ds$. By setting
\[
Y=R\wedge_\ss X, \ M=R \text{ or } Y=R, \ M=R\wedge_\ss X
\]
respectively, we obtain the spectral sequences which we call the
{\em Higher Bockstein spectral sequences}. The following theorem
summarizes the above discussion for these cases. Before we state
it, we give some comments.

The action of the Bockstein operators $Q_j\in L^*_R(L)$
(Proposition \ref{essential}) on $F_*(X)$ and on $F^*(X)$ referred
to in the statement is obtained from the natural actions of
$F^*(F)$ via the canonical algebra maps $L^*_R(L) \to F^*_R(F) \to
F^*(F)$.

The $F_*^R(L)$--comodule structure on $F_*(X)$ referred to is the
one induced by viewing it as $F_*(X)\cong F_*^R(R\wedge_\ss X)$.
If $F_*(F)$ is $F_*$--flat, it is given by composing the coaction
\[
F_*(X) \lra F_*(F)\otimes F_*(X)
\]
with the map $F_*(F)\to F_*^R(L)$ induced by
\[
(T\wedge_R L) \wedge_\ss (T\wedge_R L) \to (T_0\wedge_R L)\wedge_R
(T_1\wedge_R L) \to T_{01}\wedge_R L\wedge_R L.
\]
The $T$--module $\wh T$ was defined in \ref{that} and identified
under certain conditions as $L^R_F(T)$ and as $L_F(T)$ in
Proposition \ref{localization}. Recall that $F_*^R(L)\cong
\Lambda_{F_*}(a_0, a_1, \ldots)$ as bialgebras, if $F$ is
commutative.

\begin{thm}
\label{hbss} Let $(R, T, F)$ be a regular triple. Then there are
conditionally convergent spectral sequences of $T_*$--modules, the
cohomological and the homological {Higher Bockstein spectral
sequences}, of the form:
\begin{align}
\label{cohomhbss} E_1^{*,*}  & =  \gr^*_I(T^*) \otimes F^*(X)
\Longrightarrow  \wh T^*(X)
\\
\label{homhbss} E^1_{*,*} & =  \gr^*_I(T_*) \otimes F_*(X)
\Longrightarrow (F^\wedge_R (R\wedge_\ss X))_*
\end{align}
For finite $X$ the target of \eqref{homhbss} can be identified
with $\wh T_*(X)$. If $T=R$ and $I$ is finitely generated, then
the targets of the spectral sequences are isomorphic to
$(L_F^R(R))^*(X)$ and $(L_F^R(R\wedge_\ss X))_*$ respectively. The
differential $d_1$ of \eqref{cohomhbss} is the
$\gr_I^*(T^*)$--linear map determined by
\[
\bigoplus_{j} \bar v_j\otimes Q_j \: F^*(X) \lra I/I^2\otimes
F^*(X).
\]
Explicitly, we have for $\bar v\otimes x\in I^s/I^{s+1} \otimes
F^*(X)$
\[
d_1(\bar v\otimes x)= \sum_j \overline{v_j v}\otimes Q_j(x).
\]
A similar description holds for the differential $d^1$ of
\eqref{homhbss}. If $F$ is commutative as an $R$--ring spectrum
and $F_*^R(T)$ is $F_*$--flat, the $E_2$--term of \eqref{homhbss}
can be expressed as
\[
E^2_{*,*} \cong \Coext^{*,*}_{F_*^R(L)}(F_*, F_*(X)).
\]
If moreover $F_*(X)$ is $F_*$--free, the $E^2$-term of
\eqref{cohomhbss} is given by
\[
E_2^{*,*} \cong \Coext^{*,*}_{F_*^R(L)}(F_*(X), F_*).
\]
\end{thm}

\begin{cor}\label{even}
If $F^*(X)$ is concentrated in even degrees, the same is true for
$\wh T^*(X)$. Furthermore, the natural map $\wh T^*(X)\to F^*(X)$
is surjective.
\end{cor}

We can now give the proof of Lemma \ref{powersofideals}.

\begin{proof}[Proof of Lemma \ref{powersofideals}]
If we cut away in an obvious sense the part of the $I$--adic tower
over $T$ to the left of $T/I^s$, we obtain a spectral sequence
converging strongly to $(T/I^s)^*(X)$. In fact, it collapses at
the stage $E_s$. Let $(F^l)$ denote the filtration of
$(T/I^s)^*(X)$ coming from the spectral sequence. We show that
$I^k F^l\subset F^{k+l}\ (\ast)$. Because $F^s$ is trivial, this
will prove the lemma. To prove statement $(\ast)$, note that the
$T_*$--action on the $E_1$--term factors through $T_*/I\cong F_*$.
Hence $I$ acts trivially on the $E_r$--term for all $r$. As the
spectral sequence collapses, it follows that $I$ acts trivially on
$F^l/F^{l+1}$. This means that multiplication by an element of $I$
augments the filtration index at least by one. By induction over
$k$, this proves $(\ast)$.
\end{proof}

\begin{prop}\label{enhatnilpotent}
Let $E=\wh{E(n)}$ and $K=K(n)$ for some prime $p$. The natural map
$(L_K(E\wedge_\ss X)) \lra (L_K^E(E\wedge_\ss X))$ is an
equivalence for any spectrum $X$. Hence the homological Bockstein
spectral sequence converges for $R=T=E$ and $F=K$ conditionally to
$(L_K(E\wedge_\ss X))_*$ for any $X$.
\end{prop}

\begin{proof}
We may identify $K^\wedge_E(E\wedge_\ss
X)=\holim_s(E/I^s\wedge_\ss X)$ with the Bousfield localization
$L_K^E(E \wedge_\ss X)$ by Theorem \ref{hbss}. So the second
statement is indeed a consequence of the first. It follows from
\cite[Thm.\@ 4.2]{greenleesmay} that the Bousfield localization
$L_K^E(M)$ of an $E$--module $M$ is isomorphic to
$M^\wedge_{I_n}=F_E(K_E(I_n), M)$. Here $K_E(I_n)$ is the smash
product
\[
K_E(I_n) = K_E(v_0) \wedge_E \cdots \wedge_E K_E(v_{n-1}),
\]
where the $K_E(v_i)$ are defined by the cofibre sequences
\[
K_E(v_i) \lra E \lra E[1/v_i]
\]
with $E[1/v_i]\cong \hocolim(E \xra{v_i} E \xra{v_i} \cdots)$. We
may write $F_E(K_E(I_n), M) \cong F_{MU}(K_{MU}(I_n), M)$, where
$K_{MU}(I_n)$ is defined as $K_E(I_n)$, but $E$ is replaced by
$MU$ everywhere. Setting $M=E\wedge_\ss X$, we find that
\[
L_K^E(E\wedge_\ss X) \cong F_{MU}(K_{MU}(I_n), E\wedge_\ss X).
\]
On the other hand, we have
\[
L_K(E\wedge_\ss X) \cong F_\ss(\Gamma_{I_n}(S) \wedge_\ss L_E(S),
L_E(E\wedge_\ss X))
\]
by \cite[Prop.\@ 7.10]{hs}, where $\Gamma_{I_n}(S)$ is the
notation from \cite{greenleesmay} for the homotopy colimit over a
certain type of sequence of finite complexes of type \nolinebreak
$n$. The right hand side is isomorphic to
\[
F_\ss(\Gamma_{I_n}(S), E\wedge_\ss X) \cong
F_{MU}(MU\wedge_\ss\Gamma_{I_n}(S), E\wedge_\ss X).
\]
Now by \cite[Prop.\@ 6.6]{greenleesmay} the $p$--localizations of
$K_{MU}(I_n)$ and $MU\wedge_\ss \Gamma_{I_n}(S)$ are isomorphic.
Thus the claim follows.
\end{proof}

\subsection{Examples}

We illustrate the Higher Bockstein spectral sequence (HBSS) by
considering some examples.

\begin{exmp}[Classifying spaces of finite groups]
Let $E=\wh{E(n)}$, $I=I_n$ and $K=K(n)$ for some $n>0$ and some
prime $p$. Assume that $G$ is a finite group and let $BG$ be its
classifying space. Ravenel has shown in \cite{ravenelgroups} that
$K_*(BG)$ is always finitely generated over $K_*$. It follows that
the HBSS for $(E, E, K)$ converges strongly.

If $G$ is abelian, it follows from the results in \cite[Section
5.4]{hkr} that $E_*(BG)$ is finitely generated free over $E_*$ and
that the natural maps (Lemma \ref{powersofideals})
\[
E_*/I^s\otimes_{E_*}E_*(BG) \to (E/I^s)_*(BG)
\]
are isomorphisms. Therefore the HBSS collapses at $E^1$. This
follows also from Lemma \ref{even}, because $K_*(BG)$ is
concentrated in even degrees \cite{hkr}.

More generally, $K_*(BG)$ is trivial in odd degrees for all $G$ in
the class of ``good groups'' as defined in \cite{hkr}. So
$E_*(BG)$ is $E_*$-free for all these groups. The same is true if
$G$ is a symmetric group \cite{hunton, stricklandesymm}.

However, the classifying space of the $p$-Sylow subgroup of
$\GL_4(\F_p)$ has non-trivial odd Morava $K(n)$-theory for $p>2$
and $n\geq 2$. This was proved by Kriz \cite{kriz} for $n=2$ and
$p=3$ and by Kriz and Lee \cite{kl} in the general case.
\end{exmp}

\begin{exmp}[Moore spectra]
Let $M_p$ be the mod $p$ Moore spectrum and let us consider the
HBSS for the triple $(MU, BP\langle n\rangle, k(n))$ for variety.
Write $B=BP\langle n\rangle$, $k=k(n)$ and $I=I_n$. Let $Q_0,
\ldots, Q_{n-1}\in k^*_B(k) \subset k^*_\ss(k)$ be the Bocksteins
associated to the sequence $p, v_1, \ldots, v_{n-1}$ and put
\[
L=MU/p\wedge_{MU}\cdots\wedge_{MU} MU/v_{n-1}.
\]
The natural map
\[
k\wedge_{S} M_p \simeq k \wedge_{MU} MU/p \lra k\wedge_{MU} L
\]
induces an injection of homology groups
\[
k^S_*(M_p) \cong k_*^{MU}(MU/p) \cong \Lambda_{k_*}(a_0) \lra
\Lambda_{k_*}(a_0, \ldots, a_{n-1}) \cong k^{MU}_*(L)
\]
where the $a_j$ and their formal exterior powers are chosen as in
Lemma \ref{bocksteinhomology}. By naturality and the cited lemma,
we have $Q_0(a_0)=1$. Therefore, the $E^2$-term of the spectral
sequence is given by $E^2_{*,*} \cong \gr^*_I(B_*)/p$. As it is
trivial in odd degrees, it collapses, and we find confirmed that
$B_*(M_p)\cong B_*/p$.

Consider now the mod $p^2$ Moore spectrum $M_{p^2}$. We have
$k_*(M_{p^2})\cong\Lambda_{k_*}(a_0)$, but this time $Q_0$ is
easily seen to act trivially. The fact that for $s>1$
\[
(k/I^s)_*(M_{p^2})\cong k_*/(I^s+(p^2))
\]
implies that $d^2(\bar v\otimes a_0) = \overline{p^2v}\otimes 1$
for any monomial $v$ in $p, v_1, \ldots, v_{n-1}$. Hence the
spectral sequence collapses at $E^3_{*,*}\cong\gr_I^*(B_*)/p^2$,
and so $B_*(M_{p^2})\cong B_*/p^2$.
\end{exmp}

\begin{exmp}[Adams--Smith--Toda spectra]
Let $E=\wh{E(n)}$, $K=K(n)$ and $I=I_n$ for an odd prime $p$. Then
$M_p$ admits a self map of degree $2p-2$, which induces
multiplication by $v_1$ in $BP$-homology (see \eg\ \cite[Thm.\@
1.5]{smith}). Define $V$ to be the cofibre. Let $Q_0, \ldots,
Q_{n-1}\in K^*_E(K)\subset K^*_\ss(K)$ be the Bocksteins
corresponding to $p, v_1, \ldots, v_{n-1}$. Arguing along similar
lines as in the previous example, we find that
$K_*(V)\cong\Lambda_{K_*}(a_0, a_1)$ with $Q_0(a_0)=1$,
$Q_1(a_1)=1$, $Q_0(a_0\wedge a_1)= a_1$ and $Q_1(a_0\wedge
a_1)=-a_0$. So we may identify the complex $(E^1_{*,*}, d^1)$ with
the Koszul complex for the elements $p, v_1\in E_*$. Hence
$E^2_{*,*}\cong \gr^*_I(E_*)/(p,v_1)$, so the spectral sequence
collapses and $E_*(V)\cong E_*/(p, v_1)$, as expected.

\end{exmp}

\Addresses\recd


\begin{thebibliography}

\itemsep 1.3pt plus 1pt
\bibitem{bakainfty}
\textbf{A Baker}, \emph{{$A\sb \infty$} structures on some spectra related
  to {M}orava {$K$}-theories}, Quart. J. Math. Oxford Ser. (2) 42 (1991)
  403--419 \MR{1135302}

\bibitem{bl}
\textbf{A Baker}, \textbf{A Lazarev}, \emph{On the {A}dams spectral
  sequence for {$R$}-modules}, \agtref1{2001}{9}{173}{199}
  \MR{1823498}

\bibitem{br}
\textbf{A Baker}, \textbf{B Richter}, \emph{{$\Gamma$}-cohomology of
  rings of numerical polynomials and {$E_\infty$}-structures on {$K$}-theory},
  \arxiv{math.AT/0304473}

\bibitem{bwart}
\textbf{A Baker}, \textbf{U W{\"u}rgler}, \emph{Liftings of formal groups and
  the {A}rtinian completion of {$v\sb n\sp {-1}{\rm BP}$}}, Math. Proc.
  Cambridge Philos. Soc. 106 (1989) 511--530 \MR{1010375}

\bibitem{bw}
\textbf{A Baker}, \textbf{U W{\"u}rgler}, \emph{Bockstein operations in
  {M}orava {$K$}-theories}, Forum Math. 3 (1991) 543--560 \MR{1129998}

\bibitem{boardmanstable}
\textbf{J\,M Boardman}, \emph{Stable operations in generalized cohomology},
  from: ``Handbook of algebraic topology'', North-Holland, Amsterdam (1995)
  585--686 \MR{1361899}

\bibitem{boss}
\textbf{J\,M Boardman}, \emph{Conditionally convergent spectral sequences},
  from: ``Homotopy invariant algebraic structures (Baltimore, MD, 1998)'',
  Contemp. Math. 239, Amer. Math. Soc., Providence, RI (1999)  49--84
  \MR{1718076}

\bibitem{bousfield}
\textbf{A\,K Bousfield}, \emph{The localization of spectra with respect to
  homology}, Topology 18 (1979) 257--281 \MR{551009}

\bibitem{chr}
\textbf{J\,D Christensen}, \emph{Ideals in triangulated categories: phantoms,
  ghosts and skeleta}, Adv. Math. 136 (1998) 284--339 \MR{1626856}

\bibitem{eilenbergmoore}
\textbf{S Eilenberg}, \textbf{J\,C Moore}, \emph{Foundations of relative
  homological algebra}, Mem. Amer. Math. Soc. No. 55 (1965) 39 \MR{0178036}

\bibitem{ekmm}
\textbf{A\,D Elmendorf}, \textbf{I Kriz}, \textbf{M\,A Mandell},
  \textbf{J\,P May}, \emph{Rings, modules, and algebras in stable homotopy
  theory}, Mathematical Surveys and Monographs 47, American Mathematical
  Society, Providence, RI (1997) \MR{1417719}

\bibitem{goersseinfty}
\textbf{P\,G Goerss}, \textbf{M\,J Hopkins}, \emph{Moduli spaces of commutative
  ring spectra}, from: ``Structured ring spectra'', London Math. Soc. Lecture
  Note Ser. 315, Cambridge Univ. Press, Cambridge (2004)  151--200
  \MR{2125040}

\bibitem{greenleesmay}
\textbf{J\,P\,C Greenlees}, \textbf{J\,P May}, \emph{Completions in algebra
  and topology}, from: ``Handbook of algebraic topology'', North-Holland,
  Amsterdam (1995)  255--276 \MR{1361892}

\bibitem{hkr}
\textbf{M\,J Hopkins}, \textbf{N\,J Kuhn}, \textbf{D\,C Ravenel},
  \emph{Generalized group characters and complex oriented cohomology theories},
  J. Amer. Math. Soc. 13 (2000) 553--594 \MR{1758754}

\bibitem{hps}
\textbf{M Hovey}, \textbf{J\,H Palmieri}, \textbf{N\,P Strickland},
  \emph{Axiomatic stable homotopy theory}, Mem. Amer. Math. Soc. 128 (1997)
  \MR{1388895}

\bibitem{hs}
\textbf{M Hovey}, \textbf{N\,P Strickland}, \emph{Morava {$K$}-theories and
  localisation}, Mem. Amer. Math. Soc. 139 (1999)  \MR{1601906}

\bibitem{hunton}
\textbf{J Hunton}, \emph{The {M}orava {$K$}-theories of wreath products},
Math. Proc. Cambridge Philos. Soc. 107 (1990) 309--318 \MR{1027783}

\bibitem{husemollermoore} \textbf{D Husemoller}, \textbf{J\,C Moore},
\emph{Differential graded homological algebra of several variables},
from: ``Symposia Mathematica, Vol.  IV (INDAM, Rome, 1968/69)'',
Academic Press, London (1970) 397--429 \MR{0310040}

\bibitem{jensen}
\textbf{C\,U Jensen}, \emph{Les foncteurs d\'eriv\'es de {$\varprojlim$} et
  leurs applications en th\'eorie des modules}, Springer-Verlag, Berlin (1972)
  \MR{0407091}

\bibitem{kriz}
\textbf{I Kriz}, \emph{Morava {$K$}-theory of classifying spaces: some
  calculations}, Topology 36 (1997) 1247--1273 \MR{1452850}

\bibitem{kl}
\textbf{I Kriz}, \textbf{K\,P Lee}, \emph{Odd-degree elements in the {M}orava
  {$K(n)$} cohomology of finite groups}, Topology Appl. 103 (2000) 229--241
  \MR{1758436}

\bibitem{maclanehomology}
\textbf{S Mac~Lane}, \emph{Homology}, Classics in Mathematics, Springer-Verlag,
  Berlin (1995) \MR{1344215}

\bibitem{maclane}
\textbf{S Mac~Lane}, \emph{Categories for the working mathematician}, Graduate
  Texts in Mathematics 5, Springer-Verlag, New York (1998) \MR{1712872}

\bibitem{matsumura}
\textbf{H Matsumura}, \emph{Commutative ring theory}, Cambridge Studies in
  Advanced Mathematics 8, Cambridge University Press, Cambridge (1989)
  \MR{1011461}

\bibitem{mayeinfty}
\textbf{J\,P May}, \emph{{$E\sb{\infty }$} ring spaces and {$E\sb{\infty
  }$} ring spectra}, with contributions by
  F Quinn, N Ray, and J Tornehave, Lecture Notes in Mathematics
 577, Springer-Verlag, Berlin (1977) \MR{0494077}

\bibitem{miller}
\textbf{H\,R Miller}, \emph{On relations between {A}dams spectral sequences,
  with an application to the stable homotopy of a {M}oore space}, J. Pure Appl.
  Algebra 20 (1981) 287--312 \MR{604321}

\bibitem{ravenelgroups}
\textbf{D\,C Ravenel}, \emph{Morava {$K$}-theories and finite groups}, from:
  ``Symposium on Algebraic Topology in honor of Jos\'e Adem (Oaxtepec, 1981)'',
  Contemp. Math. 12, Amer. Math. Soc., Providence, R.I. (1982)  289--292
  \MR{676336}

\bibitem{ravenelgreen}
\textbf{D\,C Ravenel}, \emph{Complex cobordism and stable homotopy groups of
  spheres}, Pure and Applied Mathematics 121, Academic Press Inc., Orlando, FL
  (1986) \MR{860042}

\bibitem{ravenelorange}
\textbf{D\,C Ravenel}, \emph{Nilpotence and periodicity in stable homotopy
  theory}, Annals of Mathematics Studies 128, Princeton University Press,
  Princeton, NJ (1992) \MR{1192553}

\bibitem{smith}
\textbf{L Smith}, \emph{On realizing complex bordism modules. {A}pplications to
  the stable homotopy of spheres}, Amer. J. Math. 92 (1970) 793--856
  \MR{0275429}

\bibitem{stricklandesymm}
\textbf{N\,P Strickland}, \emph{Morava {$E$}-theory of symmetric groups},
  Topology 37 (1998) 757--779 \MR{1607736}

\bibitem{strickland}
\textbf{N\,P Strickland}, \emph{Products on {${\rm MU}$}-modules}, Trans. Amer.
  Math. Soc. 351 (1999) 2569--2606 \MR{1641115}

\bibitem{weibel}
\textbf{C\,A Weibel}, \emph{An introduction to homological algebra}, Cambridge
  Studies in Advanced Mathematics 38, Cambridge University Press, Cambridge
  (1994) \MR{1269324}

\bibitem{swmult}
\textbf{S W{\"u}thrich}, \emph{Multiplicative structures on $I$--adic
  towers}, in preparation

\bibitem{swshort}
\textbf{S W{\"u}thrich}, \emph{Homology of $I$--adic towers},
  \arxiv{math.AC/0411187}

\bibitem{swglasgow}
\textbf{S W{\"u}thrich}, \emph{Homology of powers of regular ideals}, Glasg.
  Math. J. 46 (2004) 571--584 \MR{2094811}

\end{thebibliography}
\end{document}